\documentclass[11pt]{article}

\usepackage{amsmath, amsfonts, amsthm, graphicx,amssymb}
\usepackage[top=21mm, bottom=22mm, left=21mm, right=21mm]{geometry}
\usepackage{comment}
\usepackage{hyperref}
\usepackage{thmtools,thm-restate}

\makeatletter
\@ifundefined{newcounteralias}{}{%
  \renewcommand\thmt@autorefsetup{%
    \@xa\def\csname\thmt@envname autorefname\@xa\endcsname
      \@xa{\thmt@thmname}%
  }%
}
\makeatother
\hypersetup{
	colorlinks=true,
	linkcolor=blue,
	filecolor=magenta,      
	urlcolor=cyan,
	citecolor=blue
}
\usepackage{authblk}
\usepackage{enumitem}
\usepackage[none]{hyphenat}
\usepackage{tikz}
\usetikzlibrary{calc, patterns,intersections, scopes}
\usetikzlibrary{positioning,arrows,shapes,decorations.markings,decorations.pathreplacing, decorations.pathmorphing,matrix,patterns}
\tikzstyle{vtx}=[circle,draw=black,fill=black,inner sep=0,minimum size=2pt,text=white,font=\footnotesize]
\usetikzlibrary{arrows.meta}
\usetikzlibrary{backgrounds}
\pgfdeclarelayer{bg}
\pgfsetlayers{bg,main}
\usepackage{algorithm}
\usepackage{algpseudocode}
\usepackage{caption}
\usepackage{dsfont}
\algrenewcommand\algorithmicrequire{\textbf{Input:}}
\algrenewcommand\algorithmicensure{\textbf{Output:}}

\newtheorem{theorem}{Theorem}[section]
\newtheorem{corollary}[theorem]{Corollary}
\newtheorem{lemma}[theorem]{Lemma}
\newtheorem{claim}[theorem]{Claim}

\usepackage{cleveref}
\crefname{lemma}{lemma}{lemmas}
\Crefname{lemma}{Lemma}{Lemmas}
\crefname{claim}{claim}{claims}
\Crefname{claim}{Claim}{Claims}
\crefname{question}{question}{questions}
\Crefname{question}{Question}{Questions}
\crefname{conjecture}{conjecture}{conjecture}
\Crefname{conjecture}{Conjecture}{Conjecture}

\theoremstyle{definition}
\newtheorem{definition}{Definition}

\newcommand{\eps}{\varepsilon}
\newcommand{\cF}{\mathcal{F}}
\newcommand{\cH}{\mathcal{H}}
\newcommand{\cT}{\mathcal{T}}

\DeclareMathOperator{\db}{\mathbf{db}}

\title{Independent sets in degree-bounded graphs}
\author{Istv\'an Tomon}
\affil{Department of Mathematics and Mathematical Statistics\\ 
	Ume\r{a} University, Ume\r{a}, Sweden}
\affil{\texttt{istvan.tomon@umu.se}}
\date{}

\begin{document}
	\sloppy
	\maketitle
	
	\begin{abstract}
		A hereditary family of graphs $\mathcal{F}$ is \emph{degree-bounded} if there exists $f:\mathbb{N}\rightarrow\mathbb{N}$ such that every $K_{t,t}$-free member of $\mathcal{F}$ has average degree at most $f(t)$. Examples include many extensively studied families such as string graphs, graphs with no induced copy of a fixed tree $T$, and graphs whose adjacency matrices have bounded factorization norm. This has led to the systematic study of this notion.  We prove two results for degree-bounded families, extending several results previously only known for particular instances.
		\begin{enumerate}
			\item If a bipartite member $G=(A,B;E)$ has at most $\frac{1}{2}|A||B|$ edges, then there exist $X\subset A$ and $Y\subset B$ with no edges between and $|X||Y|\geq (|A||B|)^{1-o(1)}$. We also construct examples with $|A|=|B|=n$ and edge density tending to 0 in which every such pair satisfies $|X||Y|=o(n^2)$.
			
			\item For every fixed $r$, every $n$-vertex $K_r$-free member of $\cF$ has chromatic number $n^{o(1)}$. This answers a question of Du and McCarty. 
		\end{enumerate}
	\end{abstract}
	
	\section{Introduction}
	
	A family $\cF$ of graphs is \emph{hereditary} if it is closed under taking induced subgraphs. A hereditary family $\cF$ is \emph{degree-bounded} if there exists a function $f:\mathbb{N}\rightarrow\mathbb{N}$ such that every $K_{t,t}$-free member $G\in \cF$ has average degree at most $f(t)$. Here, $K_{t,t}$ denotes the complete bipartite graph with vertex classes of size $t$, and $K_{t,t}$-free means that the graph does not contain $K_{t,t}$ as a (not necessarily induced) subgraph. 
	
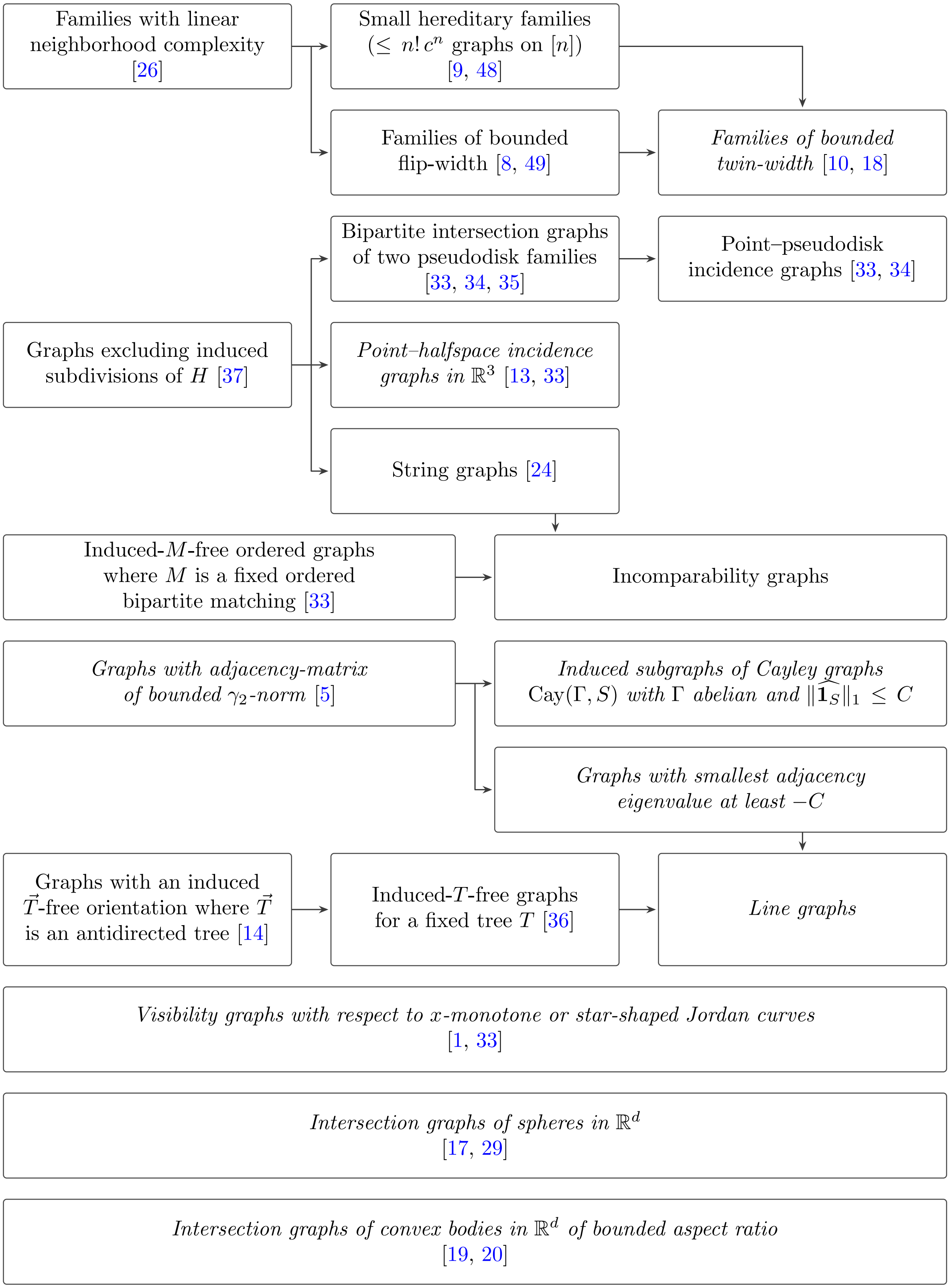
\begin{figure}\centering
\begin{tikzpicture}[
  x=1cm, y=1.75cm,
  dbfamily/.style={
    draw=black!72, line width=.48pt, rounded corners=1.5pt,
    align=center, font=\normalfont\small,
    inner xsep=4pt, inner ysep=3pt,
    minimum height=1.40cm,
    anchor=west, outer sep=0pt
  },
  dbwidth/.style={
    minimum width=#1cm,
    text width=\dimexpr#1cm-8pt\relax
  },
  dbthree/.style={dbfamily,dbwidth=\dbColWidth},
  dbtwo/.style={dbfamily,dbwidth=\dbWideWidth},
  dbfull/.style={dbfamily,dbwidth=\dbFullWidth},
  dbtall/.style={
    dbthree,minimum height=1.85cm,yshift=-.225cm
  },
  dbbottom/.style={dbfull,yshift=-.45cm},
  dblink/.style={
    draw=black!77, line width=.58pt, line join=round
  },
  dbarrow/.style={
    dblink,-{Stealth[length=4.5pt,width=3.5pt]},
    shorten >=1.2pt
  }
]
  \def\dbColWidth{4.75}
  \def\dbGap{0.65}

  \pgfmathsetmacro{\dbMidX}{\dbColWidth+\dbGap}
  \pgfmathsetmacro{\dbRightX}{2*(\dbColWidth+\dbGap)}
  \pgfmathsetmacro{\dbFullWidth}{3*\dbColWidth+2*\dbGap}
  \pgfmathsetmacro{\dbWideWidth}{(\dbFullWidth-\dbGap)/2}
  \pgfmathsetmacro{\dbWideX}{\dbWideWidth+\dbGap}
  \pgfmathsetmacro{\dbForkX}{\dbColWidth+\dbGap/2}
  \pgfmathsetmacro{\dbWideForkX}{\dbWideWidth+\dbGap/2}

  \node[dbthree] (neighborhoods) at (0,0)
    {Families with linear\\
     neighborhood complexity\\
     \cite{FPSSZ}};

  \node[dbthree] (smallclasses) at (\dbMidX,0)
    {Small hereditary families\\
     ($\le n!\,c^n$ graphs on $[n]$)\\
     \cite{BDSZ,SZZ}};

  \node[dbthree] (flipwidth) at (\dbMidX,-1)
    {Families of bounded\\
     flip-width \cite{BG,Tor}};

  \node[dbthree] (twinwidth) at (\dbRightX,-1)
    {\emph{Families of bounded}\\
     \emph{twin-width} \cite{BGKTW,DGJMR}};

  \coordinate (complexityfork) at (\dbForkX,0);
  \draw[dblink] (neighborhoods.east) -- (complexityfork);
  \draw[dbarrow] (complexityfork) -- (smallclasses.west);
  \draw[dbarrow] (complexityfork) |- (flipwidth.west);
  \draw[dbarrow] (smallclasses.east) -| (twinwidth.north);
  \draw[dbarrow] (flipwidth.east) -- (twinwidth.west);

  \node[dbthree] (pseudodisks) at (\dbMidX,-2)
    {Bipartite intersection graphs\\
     of two pseudodisk families\\
     \cite{HMST,KS,Keszegh}};

  \node[dbthree] (incidence) at (\dbRightX,-2)
    {Point--pseudodisk\\
     incidence graphs \cite{HMST,KS}};

  \node[dbthree] (subdiv) at (0,-3)
    {Graphs excluding induced\\
     subdivisions of $H$ \cite{KO}};

  \node[dbthree] (halfspaces) at (\dbMidX,-3)
    {\emph{Point--halfspace incidence}\\
     \emph{graphs in $\mathbb R^3$} \cite{CH,HMST}};

  \node[dbthree] (strings) at (\dbMidX,-4)
    {String graphs \cite{FP}};

  \coordinate (subdivfork) at (\dbForkX,-3);
  \draw[dblink] (subdiv.east) -- (subdivfork);
  \draw[dbarrow] (subdivfork) -- (halfspaces.west);
  \draw[dbarrow] (subdivfork) |- (pseudodisks.west);
  \draw[dbarrow] (subdivfork) |- (strings.west);
  \draw[dbarrow] (pseudodisks.east) -- (incidence.west);

  \node[dbtwo] (ordered) at (0,-5)
    {Induced-$M$-free ordered graphs\\
     where $M$ is a fixed ordered\\
     bipartite matching \cite{HMST}};

  \node[dbtwo] (incomp) at (\dbWideX,-5)
    {Incomparability graphs};

  \coordinate (incompentry) at ({\dbWideX+1},0);
  \draw[dbarrow]
    (strings.south -| incompentry) --
    (incomp.north -| incompentry);
  \draw[dbarrow] (ordered.east) -- (incomp.west);

  \node[dbtwo] (norm) at (0,-6)
    {\emph{Graphs with adjacency-matrix}\\
     \emph{of bounded $\gamma_2$-norm} \cite{BHT}};

  \node[dbtwo] (cayley) at (\dbWideX,-6)
    {\emph{Induced subgraphs of Cayley graphs}\\
     \emph{$\mathrm{Cay}(\Gamma,S)$ with $\Gamma$ abelian}
     \emph{and $\|\widehat{\mathbf{1}_S}\|_1\le C$}};

  \node[dbtwo] (eigenvalue) at (\dbWideX,-7)
    {\emph{Graphs with smallest adjacency}\\
     \emph{eigenvalue at least $-C$}};

  \coordinate (normfork) at (\dbWideForkX,-6);
  \draw[dblink] (norm.east) -- (normfork);
  \draw[dbarrow] (normfork) -- (cayley.west);
  \draw[dbarrow] (normfork) |- (eigenvalue.west);

  \node[dbtall] (orienttrees) at (0,-8)
    {Graphs with an induced  $\vec T$-free orientation 
     where $\vec T$ is an antidirected tree \cite{CKP}};

  \node[dbtall] (trees) at (\dbMidX,-8)
    {Induced-$T$-free graphs\\
     for a fixed tree $T$ \cite{KP}};

  \node[dbtall] (linegraphs) at (\dbRightX,-8)
    {\emph{Line graphs}};

  \draw[dbarrow] (orienttrees.east) -- (trees.west);
  \draw[dbarrow] (trees.east) -- (linegraphs.west);
  \draw[dbarrow]
    (eigenvalue.south -| linegraphs.north) --
    (linegraphs.north);

  \node[dbbottom] (visibility) at (0,-9)
    {\emph{Visibility graphs with respect to
     $x$-monotone or star-shaped Jordan curves}\\
     \cite{AK,HMST}};

  \node[dbbottom] (spheres) at (0,-10)
    {\emph{Intersection graphs of spheres in $\mathbb R^d$}\\
     \cite{DGHM,FT}};

  \node[dbbottom] (convexbodies) at (0,-11)
    {\emph{Intersection graphs of convex bodies in $\mathbb R^d$
     of bounded aspect ratio}\\
     \cite{DM25,DMN}};
\end{tikzpicture}

   \caption{A collection of known degree-bounded families, where arrows denote subfamily relation. The families in italics have linear degree-bounding functions.}
   \label{fig:deg}
 \end{figure}

	Degree-boundedness is an analogue of  another core concept in structural graph theory, called \emph{$\chi$-boundedness}. A hereditary family $\cF$ is $\chi$-bounded if there is a function $f:\mathbb{N}\rightarrow\mathbb{N}$ such that every $K_r$-free member of $\cF$ has chromatic number at most $f(r)$. This notion was introduced in an influential paper of Gy\'arf\'as \cite{Gy2} in 1987, and has been extensively studied since then, see the survey \cite{SS} as a general reference. Establishing the $\chi$-boundedness of families is generally a very difficult task. For example, the celebrated Gy\'arf\'as-Sumner conjecture \cite{Gy1,Sumner} asks whether the family of graphs avoiding a fixed tree $T$ is $\chi$-bounded, and this conjecture is still unresolved, see \cite{BB,NSS} for recent developments. However, deciding whether a family is degree-bounded is generally more manageable (while still far from trivial), and many important families of graphs have been established as degree-bounded. Interestingly, such families arise in many different contexts. Some of the main examples that are worth keeping in mind while reading this paper are the families of graphs avoiding a fixed tree $T$ \cite{KP} or all subdivisions of a fixed graph $H$ \cite{KO} as an induced subgraph, intersection and incidence graphs of certain geometric objects \cite{HMST}, and families of graphs whose adjacency matrices have bounded factorization norms \cite{BHT}. See \Cref{fig:deg} for a detailed overview of degree-bounded families. In recent years, these examples have motivated a systematic study of the structural consequences of degree-boundedness, leading to several results that reveal properties shared by all degree-bounded families.

 A result of McCarty \cite{MC21} together with a result of Kwan, Letzter, Sudakov, and Tran \cite{KLST} implies that if every $C_4$-free member of a hereditary family $\cF$ has bounded average degree, then the family is degree-bounded. Thus, in order to show that a family is degree-bounded, it is enough to consider its $K_{2,2}$-free members. The natural analogue of this for $\chi$-boundedness fails: there exist families $\cF$ in which every $K_3$-free member has bounded chromatic number, while the chromatic number of $K_4$-free members of $\cF$ may be arbitrarily large \cite{CHMS}. Another strong property of degree-boundedness is established by Gir\~ao and Hunter \cite{GH}, who proved that if $\cF$ is degree-bounded, then the function $f$ can be chosen to be a polynomial. Finally, Hunter, Milojevi\'c, Sudakov, and Tomon \cite{HMST} proved that the function $f$ can be chosen to be linear if and only if the family satisfies the \emph{density Erd\H{o}s-Hajnal property}, i.e. every member of positive edge density contains a linear sized complete bipartite graph. For an even more detailed overview of degree-bounded families, we refer the reader to the survey of Du and McCarty \cite{DM25}.

    The goal of our paper is to prove that degree-boundedness ensures even more structure. We show that such families are almost $\chi$-bounded, meaning that their $K_r$-free members can be properly colored with a subpolynomial number of colors. Moreover, bipartite degree-bounded graphs contain almost linear sized anticomplete pairs across the vertex classes.   Given a graph $G$, let $\db(G)$ denote the maximum average degree of a $C_4$-free induced subgraph of $G$. By the characterization above, a hereditary family $\cF$ is degree-bounded if and only if there is a constant $k$ such that $\db(G)\leq k$ for every $G\in \cF$. We therefore formulate our results for individual graphs with $\db(G)\leq k$.
	
	\subsection{Large anticomplete pairs}
	
	A central problem of structural graph theory is the Erd\H{o}s-Hajnal conjecture \cite{EH77,EH89}, which states that  if $\cF$ is a \emph{proper} hereditary family of graphs (meaning not the family of all graphs), then there exists $c>0$ such that every $G\in \cF$ contains either a clique or an independent set of size at least $|V(G)|^{c}$. This conjecture is still wide open, see \cite{BNSS} for recent developments. 
    
    A natural relaxation asks for two large disjoint vertex sets $X$ and $Y$ such that either all or none of the edges between $X$ and $Y$ are present. We call $(X,Y)$ a \emph{complete pair} if all edges between $X$ and $Y$ are present, \emph{anticomplete pair} if there are no edges between, and \emph{homogeneous} if it is either complete or anticomplete. The size of the pair is $\min\{|X|,|Y|\}$. Erd\H{o}s, Hajnal and Pach \cite{EHP} proved the analogue of the Erd\H{o}s-Hajnal conjecture that for every proper hereditary family $\cF$, there exists $c>0$ such that every $G\in \cF$ contains a homogeneous pair of size at least $|V(G)|^c$. This exponent $c$ can be arbitrarily small, so it is an intriguing problem to identify which families $\cF$ guarantee  homogeneous pairs of linear size. A hereditary family has the \emph{strong Erdős–Hajnal property} if every $n$-vertex member contains homogeneous pairs of size $\Omega(n)$. The name reflects the fact that this property implies the Erd\H{o}s-Hajnal property \cite{APPRS}. For results on the strong Erd\H{o}s-Hajnal property, see \cite{APPRS,CSSS1,CSSS2,FPT1,FPT2}.

Even among degree-bounded families, the exponent in the theorem of
Erd\H{o}s, Hajnal and Pach can be arbitrarily small. To see this, fix
$0<\gamma<1$ and consider the random graph $G_{n,p}$ with
$p=1-n^{-\gamma}$. With high probability, every homogeneous pair has
size $O(n^\gamma\log n)$, and every induced tree has $O(1/\gamma)$
vertices. Thus, for a sufficiently large fixed tree $T$, these graphs
belong with high probability to the family of induced $T$-free graphs,
which is degree-bounded by a theorem of Kierstead and Penrice \cite{KP}.
These examples are very dense. In contrast, for sparse graphs, Chudnovsky, Scott,
Seymour and Spirkl \cite{CSSS1} proved that for every tree $T$ there
exists $\tau>0$ such that every $n$-vertex induced $T$-free graph
of edge density at most $\tau$ contains an anticomplete pair of size at
least $\Omega(n)$. The same authors
proved an analogous result for graphs containing no induced subdivision
of a fixed graph $H$ \cite{CSSS2}.

    For general degree-bounded families, even vanishing edge density is not enough to ensure linear sized anticomplete pairs. In \Cref{lemma:simple_construction}, for every $c>0$ we construct $n$-vertex graphs $G$ with $O(n^{1+c})$ edges, $\db(G)\leq 6$, and such that every anticomplete pair has size at most $n^{1-\Omega(c)}$. While it might be tempting to concede at this stage and conclude that degree-boundedness does not force large homogeneous pairs, the situation can be salvaged by considering a bipartite setting.  Given a bipartite graph
$G=(A,B;E)$, we seek an anticomplete pair $(X,Y)$ with $X\subseteq A$
and $Y\subseteq B$. Scott, Seymour and Spirkl \cite{SSS1} proved that if $G=(A,B;E)$ is a bipartite graph which contains no induced copy of a tree $T$, and $|E|\leq (1-\gamma)|A||B|$ for some $\gamma>0$, then there is an anticomplete pair $(X,Y)$ with $X\subset A$, $Y\subset B$, and $|X||Y|\geq \Omega(|A||B|)$. For further results on homogeneous pairs in bipartite graphs, see \cite{ATW,BHT,FPS1,HMST,SSS2}. Our first theorem shows that degree-boundedness guarantees an
almost linear bound in this setting.
	
	\begin{theorem}\label{thm:main_EH}
		Fix an integer $k\geq 2$ and $\gamma>0$. Let $G=(A,B;E)$ be a bipartite graph such that $|E|\leq (1-\gamma)|A||B|$ and $\db(G)\leq k$. Then there exists an anticomplete pair $(X,Y)$ with $X\subset A$ and $Y\subset B$ such that $$|X||Y|\geq (|A||B|)^{1-o(1)},$$ where $o(1)\rightarrow 0$ as $|A||B|\rightarrow\infty$.
	\end{theorem}
	
	The almost linear bound in \Cref{thm:main_EH} cannot be strengthened
to a linear one, even for balanced bipartite graphs of vanishing density.
	
	\begin{restatable}{theorem}{construction}\label{thm:construction}
		For every $c>0$ and infinitely many $n$, there exists a bipartite graph $G=(A,B;E)$ such that $|A|=|B|=n$, $\db(G)\leq 28$, $|E|\leq n^{1+c}$, and for every anticomplete pair $(X,Y)$ in $G$ with $X\subset A, Y\subset B$, we have $$|X||Y|=O_c\left(\frac{n^2}{\log \log n}\right).$$
	\end{restatable}

    The proof of \Cref{thm:main_EH} is presented in \Cref{sect:EH}, and the construction proving \Cref{thm:construction} is presented in \Cref{sect:constr}.
	
	\subsection{Subpolynomial chromatic number}

    Although degree-boundedness and $\chi$-boundedness are analogous, neither
implies the other. For one direction, the family of all bipartite graphs is $\chi$-bounded, yet it contains $C_4$-free members of arbitrarily large average degree. The other direction is more subtle. Recall the result of K\"uhn and Osthus \cite{KO} stating that the family of graphs avoiding all subdivisions of a graph $H$ as an induced subgraph is degree-bounded. But if $H$ is the 1-subdivision of $K_5$, a construction of Pawlik, Kozik, Krawczyk, Laso\'n, Micek, Trotter, and Walczak \cite{Pawlik} shows that  this family is not $\chi$-bounded.

 Despite this, degree-boundedness still has a strong effect on the chromatic numbers of $K_r$-free members. A natural relaxation of $\chi$-boundedness is to consider whether every $n$-vertex $K_r$-free member of a hereditary family $\cF$ has chromatic number $n^{o(1)}$ for fixed $r$. This property has been  established for several special families of degree-bounded graphs such as intersection graphs of segments \cite{Suk}, string graphs \cite{FP}, families of graphs avoiding a fixed tree $T$ as an induced subgraph \cite{BB,NSS}, families of graphs avoiding all induced subdivisions of a fixed graph $H$ \cite{NSS2}. Du and McCarty \cite[Problem 3.9]{DM25} asked whether degree-boundedness already forces subpolynomial chromatic number in triangle-free graphs.  Our next theorem proves the stronger result that a subpolynomial bound holds in $K_r$-free members for every fixed $r$.
    
	\begin{theorem}\label{thm:main_color}
		Fix positive integers $r$ and $k$. If $G$ is an $n$-vertex $K_r$-free graph such that $\db(G)\leq k$, then $\chi(G)=n^{o(1)}$.
	\end{theorem}

    The proof of \Cref{thm:main_color} is presented in \Cref{sect:chromatic}.

    \section{Proof overview}

    We outline the proofs of \Cref{thm:main_EH} and \Cref{thm:main_color}. As the two proofs are very similar, we focus mainly on the proof of \Cref{thm:main_EH}.

\medskip

\noindent
   \textbf{Outline of the proof of \Cref{thm:main_EH}.} Let $G=(A,B;E)$ be a bipartite graph such that $\db(G)\leq k$ and $|E|\leq (1-\gamma)|A||B|$. Then we would like to find an anticomplete pair $(X,Y)$ with $X\subset A$, $Y\subset B$ and $|X||Y|\geq (|A||B|)^{1-o(1)}$. If $G$ contains no induced $C_4$-free subgraph with large average degree, then we expect most pairs of neighborhoods of vertices of $B$ to either heavily overlap, or have little overlap. This property can be captured by constructing a sparse hypergraph $\cH$ on $A$ such that every neighborhood of a vertex of $B$ is almost completely covered by edges of $\cH$. Then we find our anticomplete pair by taking subsets of $A$ and $B$ successively so that $\cH$ is sparsified in each round.   We will show that it is enough to consider balanced bipartite graphs, so we assume $|A|=|B|=n$. Moreover, by a VC-dimension argument, we may further assume that the maximum degree of $G$ is at most $n^{1-\tau}$.
    
    Formally, the proof hinges on the following definition. Say that $G$ has the \emph{$(L,D,s)$-property} if there is a mixed hypergraph $\cH$ on $A$ with edges of sizes $h=2,\dots,s$ such that the number of edges of size $h$ is
    $$|\cH_h|\leq LnD^{h-1},$$
    and $\cH[N(x)]$ has no independent set of size $s$ for every $x\in B$. Then the proof is the combination of the following two auxiliary lemmas.
    \begin{enumerate}
        \item If $\db(G)\leq k$, then $G$ has the $(O(1),\Delta(G),k)$-property (\Cref{lemma:k-property}).

        \textbf{Proof outline.} We define a fractional version of the requirements, which leads to a linear program. We prove that if the optimal solution of this program is sufficiently small, then rounding gives the desired hypergraph. On the other hand, if the optimal solution is too large, then an optimal solution of the dual leads to a $C_4$-free induced subgraph with large average degree.
        
        \item If a bipartite graph $G=(A,B;E)$ has the $(L,D,s+1)$-property with $D=\Delta(G)$ and is sufficiently regular, then it has a large induced subgraph with the $(L',D,s)$-property, where $L'$ is bounded by some function of $L$ and the measure of regularity (\Cref{lemma:s-property}).

        \textbf{Proof outline.} We consider the graph $\cH_2$ on the vertex set $A$. Sample random subsets $S^*\subset S\subset A$. We keep those vertices of $A\setminus S$ that have no $\cH_2$-neighbor in $S^*$, and we keep those vertices of $B$ that have a unique $G$-neighbor in $S$, and this neighbor is in $S^*$.  We argue that by setting the sampling probabilities optimally, this gives the desired induced subgraph with positive probability. 
    \end{enumerate} 
    Thus, by the first lemma, we start with a graph having the $(O(1),\Delta(G),k)$-property. Then, we iterate the second lemma $k-1$ times until we arrive at a large induced subgraph having the $(L,D,1)$-property. But a bipartite graph with the $(L,D,1)$-property is empty, giving the desired anticomplete pair. One subtle issue with this argument is that the second lemma requires the input graph to be sufficiently regular, which might not be inherited by the output graph. Therefore, we also combine this procedure with a density decrement argument, and we argue that either regularity is preserved, or we can pass to a much sparser subgraph.

    \medskip

    \noindent
    \textbf{Outline of the proof of \Cref{thm:main_color}.} Let $G$ be an $n$-vertex $K_r$-free graph such that $\db(G)\leq k$. We want to prove that $\chi(G)=n^{o(1)}$. We proceed by induction on $r$, so observing that the neighborhood of every vertex is $K_{r-1}$-free, we may assume that the neighborhood of every vertex of $G$ is already $K$-colorable with $K=n^{o(1)}$. Our goal is to show that $G$ contains an independent set of size $n^{1-o(1)}$, as then repeatedly taking out large independent sets leads to a proper coloring of $G$ with few colors. We introduce the following variant of the $(L,D,s)$-property.

     The graph $G$ has the \emph{stable-$(L,D,s)$-property} if there is a mixed hypergraph $\cH$ on $V(G)$ with edges of sizes $h=2,\dots,s$ such that the number of edges of size $h$ is
    $$|\cH_h|\leq LnD^{h-1},$$
    $\cH[N(x)]$ has no independent set of size $s$ for every $x\in V(G)$, and every edge of $\cH$ is an independent set in $G$. Then the proof is the combination of the following two auxiliary lemmas, similarly as above.
    \begin{enumerate}
        \item If $\db(G)\leq k$, $G$ is $K_r$-free, and the neighborhood of every vertex is $K$-colorable, then $G$ has a large induced subgraph with the stable-$(O_k(K^k),\Delta(G),O_{k,r}(1))$-property (\Cref{lemma:stable}).
        \item If  $G$ has the stable-$(L,D,s+1)$-property and is sufficiently regular, then it has a large induced subgraph with the stable-$(L',D,s)$-property, where $L'$ is bounded by a function of $L$ and the measure of regularity (\Cref{lemma:stable-s-decrease}).
    \end{enumerate} 
    Thus, by the first lemma, we start with a graph having the stable-$(n^{o(1)},\Delta(G),O_{r,k}(1))$-property. Then, we iterate the second lemma a constant number of times until we arrive at a large induced subgraph having the stable-$(n^{o(1)},D,1)$-property. This property implies that the graph is empty, giving the desired independent set. Keeping track of the regularity required for the second lemma will cause issues, which we resolve by a density decrement argument.

	\section{Preliminaries}

    Here, we present our notation, which is mostly standard, and collect a few basic results needed for our proofs.

    \medskip

    \noindent
    \textbf{Graph theoretic notation.}  If $G$ is a graph and $X\subset V(G)$, then $G[X]$ is the subgraph of $G$ induced on $X$. We write $e(G)$ for the number of edges of $G$, $d(G)$ for its average degree, and $\Delta(G)$ for its maximum degree. Given a vertex $x\in V(G)$, $N_G(x)$ is the neighborhood of $x$, or we simply write $N(x)$ if $G$ is clear from the context. We write $\chi(G)$ for the chromatic number of $G$, $\omega(G)$ for its clique number, and $\alpha(G)$ for its independence number. If $G=(A,B;E)$ is a bipartite graph and $X\subset A$, $Y\subset B$, we write $G[X,Y]$ instead of $G[X\cup Y]$, and $e(X,Y)$ instead of $e(G[X,Y])$, if $G$ is clear from the context. 
	
	We also use the following notation regarding (mixed) hypergraphs. Let $\cH$ be a hypergraph, and let $V(\cH)$ denote the set of vertices and $E(\cH)\subset 2^{V(\cH)}$ the set of edges. Given $X\subset V(\cH)$, $\cH[X]$ is the subhypergraph of $\cH$ induced on $X$, that is, we keep only those edges of $\cH$ that are completely contained in $X$. For a positive integer $h$, we write $\cH_h$ for the set of edges of size exactly $h$. The independence number of $\cH$, denoted by $\alpha(\cH)$, is the size of the largest set $X$ such that $\cH[X]$ is empty.

\medskip

\noindent
\textbf{Linear programming.} The proofs of our main results rely on linear programming duality. We consider the primal program
    $$\text{Minimize } \mathbf{b}^T \mathbf{x}\quad\text{subject to }A\mathbf{x}\geq \mathbf{c}, \mathbf{x}\geq 0.$$
    Its \emph{dual linear program} is 
    $$\text{Maximize }\mathbf{c}^T \mathbf{y}\quad\text{subject to }A^T\mathbf{y}\leq \mathbf{b}, \mathbf{y}\geq 0.$$
    The \emph{strong duality theorem} states that if one of the problems has an optimal solution, then so does the other, and $\min_{\mathbf{x}} \mathbf{b}^T\mathbf{x}=\max_{\mathbf{y}} \mathbf{c}^T\mathbf{y}$. 

    \medskip

    \noindent
    \textbf{A simple probabilistic lemma.}
    
	\begin{lemma}\label{lemma:simple_prob}
		Let $M>0$, and let $Z\geq 0$ be a random variable such that $Z\leq M$. Then for $c\in (0,1)$,
		$$\mathbb{P}(Z\geq c\mathbb{E}Z)\geq \frac{(1-c)\mathbb{E}Z}{M}.$$
	\end{lemma}
	
	\begin{proof}
		Let $t=c\mathbb{E}Z$, then
		$$\mathbb{E}Z\leq M\mathbb{P}(Z\geq t)+t\mathbb{P}(Z<t)\leq M\mathbb{P}(Z\geq t)+t,$$
		from which
		$$\mathbb{P}(Z\geq t)\geq \frac{\mathbb{E}Z-t}{M}=\frac{(1-c)\mathbb{E}Z}{M}.$$
	\end{proof}

    \noindent
    \textbf{Degree-bounded graphs with no large anticomplete pairs.} As promised in the introduction, we show the existence of sparse degree-bounded graphs with no linear or almost linear anticomplete pairs. As the proof is based on standard probabilistic arguments, we only give a sketch.

    \begin{lemma}\label{lemma:simple_construction}
    For every $c\in (0,1)$, if $n$ is sufficiently large, there exists an $n$-vertex graph $G$ such that $e(G)\leq O(n^{1+c})$, $\db(G)\leq 6$, and every anticomplete pair in $G$ has size at most $O(n^{1-c/3}\log n)$.
    \end{lemma}

    \begin{proof}
    Let $G_1$ be the union of $n^{1-c}$ cliques of size $n^c$, and let $G_2$ be a copy of the Erd\H{o}s-R\'enyi random graph $G_{n,p}$ with $p=\frac{1}{4}n^{-1+c/3}$ on the same vertex set as $G_1$. Let $G$ be the union of $G_1$ and $G_2$. With high probability, $e(G)=O(n^{1+c})$. Note that if $X\subset V(G)$ induces a $C_4$-free graph, then $X$ intersects every clique of $G_1$ in at most 3 vertices. Therefore, $|X|\leq 3n^{1-c}$. Next, we show that every such $X$ induces at most $2|X|$ edges in $G_2$ with high probability. Indeed, for $m\leq 3n^{1-c}$, we have
    $$\mathbb{P}(\exists X : |X|=m, e(G_2[X])\geq 2m)\leq \binom{n}{m}\binom{\binom{m}{2}}{2m}p^{2m}\leq \left(\frac{en}{m}\right)^m \left(\frac{em}{4}\right)^{2m}p^{2m}=\left(\frac{e^3p^2nm}{16}\right)^m<n^{-cm/3}.$$
    Thus, with high probability, for every $X\subset V(G)$ such that $G[X]$ is $C_4$-free, we have $$e(G[X])\leq e(G_1[X])+e(G_2[X])\leq 3|X|.$$ This shows that the average degree of $G[X]$ is at most $6$. Finally, again from standard concentration arguments, every anticomplete pair in $G_2$ has size $O(n^{1-c/3}\log n)$ with high probability, thus the same holds for $G$.
    \end{proof}
	
	\section{Degree-bounded bipartite graphs without linear anticomplete pairs}\label{sect:constr}
	
	In this section, we prove \Cref{thm:construction}, which we recall here for the reader's convenience.

    \construction*
    
    Our graph is constructed by taking the union of the blow-ups of random matchings at different scales. Let $r$ be a positive integer, and define the sequence $t_1,\dots,t_r$ as follows. Let $t_1=100r^2$, and for $i=1,2,\dots,r-1$, define $t_{i+1}=t_i^2$. Let $n=t_{r}^{\lceil 2/c\rceil}$, and write $n_i=n/t_i$. We observe that
	$$r=\Theta_c(\log \log n).$$
	
	Let $G$ be a bipartite graph, whose vertex classes, $A$ and $B$, are identified with $[n]$. For $i=1,\dots,r$, let $(A_{i,j})_{1\leq j\leq n_i}$ be the partition of $A$ into intervals of size $t_i$, and let $(B_{i,j})_{1\leq j\leq n_i}$ be a partition of $B$ into intervals of size $t_i$. We call these sets  \emph{level-$i$ intervals}.  Let $\pi_i:[n_i]\rightarrow [n_i]$ be a permutation chosen uniformly at random, independent of all other choices, and let $G$ be the graph in which we add a complete bipartite graph between $A_{i,j}$ and $B_{i,\pi_i(j)}$ for every $i\in [r]$ and $j\in [n_i]$. We call the complete bipartite graph between  $A_{i,j}$ and $B_{i,\pi_i(j)}$ a \emph{level-$i$ block}. 
	
	We prove that with positive probability, $G$ satisfies the required conditions. First, note that the number of edges of $G$ is at most
	$$\sum_{i=1}^r nt_i\leq 2nt_r\leq n^{1+c}.$$
	Next, we consider the average degrees of $C_4$-free induced subgraphs of $G$.
	
	\begin{lemma}\label{lemma:construction_C4}
		With probability at least $3/4$, we have $\db(G)\leq 28$.
	\end{lemma}
	
	\begin{proof}
		Let $S\subset A$ and $T\subset B$, then we consider the probability of the event that $G[S,T]$ is $C_4$-free and has average degree more than 28. Say that an edge $e$ between $S$ and $T$ is \emph{lonely} if $e$ connects some intervals $A_{i,j}$ and $B_{i,\pi_i(j)}$ satisfying $|A_{i,j}\cap S|=1$ and $|B_{i,\pi_i(j)}\cap T|=1$. First, we show that only a few edges of $G[S,T]$ are not lonely if $G[S,T]$ is $C_4$-free,
		
		\begin{claim}\label{claim:lonely}
			If $G[S,T]$ is $C_4$-free, then the number of edges of $G[S,T]$ that are not lonely is at most $|S|+|T|$.
		\end{claim}
		
		\begin{proof}
			For every $u\in S$ and $i\in [r]$, let $V_{u,i}$ be the set of vertices $v\in T$ such that the edge $uv$ is in a level-$i$ block, and not in a level $j$-block for $j>i$. The main observation is that the family $\{V_{u,i}:u\in S,i\in [r]\}$ is  \emph{laminar}, that is, any two of its members are either disjoint, or one contains the other. But then $C_4$-freeness implies that the sets $V_{u,i}$, where $|V_{u,i}|\geq 2$, are pairwise disjoint. Thus,
			$$\sum_{u\in S,i\in [r], |V_{u,i}|\geq 2} |V_{u,i}|\leq |T|.$$
			Similarly, for $v\in T$ and $i\in [r]$, let $U_{v,i}$ be the set of vertices $u\in S$ such that the edge $uv$ is in a level-$i$ block, and not in a level $j$-block for $j>i$. Then, by the same reasoning,
			$$\sum_{v\in T,i\in [r], |U_{v,i}|\geq 2} |U_{v,i}|\leq |S|.$$
			But every edge that is not lonely is of the form $uv$ for some $v\in V_{u,i}$, $|V_{u,i}|\geq 2$, or $u\in U_{v,i}$, $|U_{v,i}|\geq 2$, so the number of edges that are not lonely is at most $|S|+|T|$.
		\end{proof}
		
		Next, our goal is to bound the number of lonely edges. Let $X_i$ denote the number of lonely edges added with a level-$i$ block, that is,
		$$X_i=|\{j\in [n_i]: |S\cap A_{i,j}|=|T\cap B_{i,\pi_i(j)}|=1\}|.$$ Let $S_i$ be the set of vertices $u\in S$ such that $\{u\}=S\cap A_{i,j}$ for some $j\in [n_i]$, and define $T_i$ analogously. Then $X_i$ is a random variable with hypergeometric distribution of parameters $(n_i,|S_i|,|T_i|)$. Hence, we can bound the moments of $X_i$ as follows:
		$$\mathbb{E}\binom{X_i}{q}=\frac{\binom{|S_i|}{q}\binom{|T_i|}{q}}{\binom{n_i}{q}}\leq \frac{(|S_i||T_i|/n_i)^q}{q!}\leq \frac{(|S||T|/n_i)^q}{q!}.$$
		For $I\subset [r]$, write $X_I=\sum_{i\in I}X_i$, and let $X=X_{[r]}$. Then $X$ is an upper bound on the number of lonely edges. Using the independence of the variables $X_i$, $i\in [r]$, we have
		$$\mathbb{E}\binom{X_I}{q}=\sum_{\sum_{i\in I}q_i=q}\prod_{i\in I}\mathbb{E}\binom{|X_i|}{q_i}\leq \sum_{\sum_{i\in I}q_i=q} \prod_{i\in I}\frac{(|S||T|/n_i)^{q_i}}{q_i!}=\frac{(|S||T|\sum_{i\in I}n_i^{-1})^q}{q!}.$$
		From this, we get the tail bound
		$$\mathbb{P}(X_I\geq q)\leq \mathbb{E}\binom{|X_I|}{q}\leq \frac{(|S||T|\sum_{i\in I}n_i^{-1})^q}{q!}.$$
		We divide the indices $i$ into two sets depending on the sizes of $n_i$. In what follows, assume that $|S|\leq |T|$.  Let $i_0=r+1$ if $n_r>|T|$, otherwise, let $i_0$ be the smallest index $i$ such that $n_i\leq |T|$. Set $J=\{i\in [r]:i\geq i_0-1\}$. By observing the trivial bound $X_i\leq \min\{|S|,n_i\}$, we can bound $X_J$ as 
		$$X_J\leq |S|+\sum_{i=i_0}^{r} n_i\leq 3|T|,$$ 
		where the sum is empty in case $i_0=r+1$. In the second inequality, we use that $n_{i_0}\leq |T|$ and $n_{i+1}\leq n_i/2$. On the other hand, setting $I=[r]\setminus J=\{i\in [r]:i\leq i_0-2\}$, the previous tail bound gives
		$$\mathbb{P}(X_I\geq q)\leq \frac{(|S||T|\sum_{i\in I}n_i^{-1})^q}{q!}\leq \left(\frac{3|S||T|\sum_{i\in I}n_i^{-1}}{q}\right)^q.$$
		If $i_0\leq 2$, the right hand side is 0. Otherwise, let $t=t_{i_0-2}$, then $|T|<n_{i_0-1}=n/t^2$ and
		$$\sum_{i\in I}n_i^{-1}\leq 2n_{i_0-2}^{-1}=\frac{2t}{n}.$$
		Hence,
		$$\mathbb{P}(X_I\geq 9|T|)\leq \left(\frac{3|S||T|\sum_{i\in I}n_i^{-1}}{9|T|}\right)^{9|T|}\leq \left(\frac{2|T|t}{3n}\right)^{9|T|}.$$
		Fix some positive integer $b\leq n$. Then the expected number of pairs $(S,T)$ with $|S|\leq |T|=b$ and $X_I\geq 9b$ is at most  
		$$\binom{n}{\leq b}\binom{n}{b}\left(\frac{2bt}{3n}\right)^{9b}\leq \left(\frac{9n^2}{b^2}\right)^b \left(\frac{2bt}{3n}\right)^{9b}\leq \left(\frac{9b^7t^9}{n^7}\right)^b\leq \left(\frac{9}{t^5}\right)^b<10^{-b}.$$
		In the third inequality, we used that $b<n/t^2$, and in the last, that $t\geq t_1>10$. By taking union bound over $b=1,\dots,n$, we get that with probability at least $8/9$, for every pair $(S,T)$ satisfying $1\leq |S|\leq |T|$, we have $X_I\leq 9|T|$, and thus 
		$$X= X_I+X_J\leq 12|T|.$$
		Hence, the number of lonely edges of $G[S,T]$ is at most $12|T|$. Combining this with \Cref{claim:lonely}, if $G[S,T]$ is $C_4$-free, then its number of edges is at most $14|T|$, showing $d(G[S,T])\leq 28$. By considering the symmetric case $|S|\geq |T|$ as well, this shows that with probability at least $3/4$, every $C_4$-free induced subgraph of $G$ has average degree at most 28, finishing the proof.
	\end{proof}
	
	Next, we consider the size of anticomplete pairs in $G$.
	
	\begin{lemma}\label{lemma:construction_biclique}
		With probability at least $3/4$, every anticomplete pair $(X,Y)$ in $G$ with $X\subset A$ and $Y\subset B$ satisfies
        $$|X||Y|\leq \frac{4n^2}{r}.$$
	\end{lemma}
	
	\begin{proof}
		First, we show that for every $U\subset A$ and $V\subset B$ that are unions of level-$i$ intervals, the number of level-$(i-1)$ blocks between $U$ and $V$ is close to its expectation. 
		
		Let $i\in \{2,\dots,r\}$, and let $U$ be a union of level-$i$ intervals in $A$, and let $V$ be a union of level-$i$ intervals in $B$. Write $\alpha=\alpha_U=|U|/n$ and $\beta=\beta_V=|V|/n$. Let $Z=Z_{U,V}$ be the number of level-$(i-1)$ blocks between $U$ and $V$, then $Z$ has hypergeometric distribution with parameters $(n_{i-1},\alpha n_{i-1},\beta n_{i-1})$. Thus, $\mathbb{E}Z=\alpha\beta n_{i-1}$. Applying Hoeffding's inequality for random samples without replacement (see e.g. \cite[Proposition 1.2]{BM}), we have for $t\geq 0$ that
		$$\mathbb{P}(Z\leq \mathbb{E}Z-t)\leq \exp\left(-\frac{2t^2}{\alpha n_{i-1}}\right)\leq \exp\left(-\frac{2t^2}{n_{i-1}}\right).$$
		Hence, by setting $t=n_{i-1}/r$, we get that
		$$\mathbb{P}(Z\leq (\alpha\beta-1/r)n_{i-1})\leq \exp\left(-\frac{2n_{i-1}}{r^2}\right).$$
		Say that the pair $(U,V)$ is \emph{slim} if $Z\leq (\alpha\beta-1/r)n_{i-1}$. The number of choices for $(U,V)$ is at most $4^{n_i}$, so by the union bound, 
		$$\mathbb{P}\left(\exists U,V: (U,V)\text{ is slim}\right)\leq 4^{n_i}\exp\left(-\frac{2n_{i-1}}{r^2}\right)\leq \exp\left(2n_i-\frac{2n_{i-1}}{r^2}\right)\leq 10^{-i}.$$
		The last inequality follows as $n_i=n_{i-1}/t_{i-1}$ and $t_{i-1}\geq 100r^2 2^{i-2}$. Thus, by taking union bound for $i=2,\dots,r$, we get that with probability at least $3/4$, for every $i\in \{2,\dots,r\}$, every pair $(U,V)$, where $U\subset A$ and $V\subset B$ is a union of level-$i$ intervals, is not slim. Let $\mathcal{E}$ denote this event. We prove that $\mathcal{E}$ ensures the required bound on the size of anticomplete pairs.
		
		Let $X\subset A$ and $Y\subset B$ such that there are no edges between $X$ and $Y$. Let $\alpha=|X|/n$ and $\beta=|Y|/n$. For $i\in [r]$, let $X_i$ and $Y_i$ be the union of level-$i$ intervals that intersect $X$ and $Y$, respectively, and let $\alpha_i=|X_i|/n$ and  $\beta_i=|Y_i|/n$. We note that $\alpha_i\geq \alpha_{i-1}\geq \alpha$ and $\beta_i\geq \beta_{i-1}\geq \beta$. 
		
		Let $i\in \{2,\dots,r\}$, and let $L$ be the number of level-$(i-1)$ blocks that intersect $X_{i-1}$ or $Y_{i-1}$. As there are no edges between $X$ and $Y$, and the blocks form a matching of the level-$(i-1)$ intervals, we have  
		$$L=(\alpha_{i-1}+\beta_{i-1})n_{i-1}.$$ Now let $M$ be the number of level-$(i-1)$ blocks that intersect $X_i$ or $Y_i$. Clearly $L\leq M$. Moreover, 
		$$M\leq \left(\alpha_i+\beta_i-\alpha_i\beta_i+\frac{1}{r}\right)n_{i-1}.$$
		Indeed, $\alpha_in_{i-1}$ and $\beta_i n_{i-1}$ are the numbers of level-$(i-1)$ intervals contained in $X_i$ and $Y_i$, respectively. Each such interval is a vertex class of a level-$(i-1)$ block, and those blocks that are contained in $G[X_i,Y_i]$ are counted twice. By the event $\mathcal{E}$, the number of such blocks is at least $(\alpha_i\beta_i-1/r)n_{i-1}$.
		
		Thus, comparing our bound on $L$ and $M$, we have
		$$\alpha_{i-1}+\beta_{i-1}\leq \alpha_i+\beta_{i}-\alpha_i\beta_i+\frac{1}{r}\leq \alpha_i+\beta_i+\left(\frac{1}{r}-\alpha\beta\right),$$
		or equivalently,
		$$\alpha\beta-\frac{1}{r}\leq \alpha_i+\beta_i-(\alpha_{i-1}+\beta_{i-1}).$$
		Summing these inequalities for $i=2,\dots,r$, we get
		$$(r-1)\alpha\beta-\frac{r-1}{r}\leq \alpha_r+\beta_r-(\alpha_1+\beta_1)\leq 2,$$
		from which we get $\alpha\beta< 4/r$.
	\end{proof}
	
	Combining \Cref{lemma:construction_C4} and \Cref{lemma:construction_biclique}, we conclude that with probability at least $1/2$, the graph $G$ satisfies the required conditions of \Cref{thm:construction}.
	
	\section{Anticomplete pairs in degree-bounded bipartite graphs}\label{sect:EH}
	
	In this section, we prove \Cref{thm:main_EH}. In  particular, we deduce it from the following theorem, which is the main result of this section.
    
    \begin{theorem}\label{thm:main_EH_balanced}
    Let $k\geq 2$ be an integer and $\gamma,\eps>0$, then the following holds if $n$ is sufficiently large. Let $G=(A,B;E)$ be a bipartite graph such that $|A|=|B|=n$, $|E|\leq (1-\gamma)n^2$, and $\db(G)\leq k$. Then there exists an anticomplete pair $(X,Y)$ such that $X\subset A$, $Y\subset B$, and $|X|,|Y|\geq n^{1-\eps}$.
    \end{theorem}

    \begin{proof}[Proof of \Cref{thm:main_EH} assuming \Cref{thm:main_EH_balanced}]
    Let $G=(A,B;E)$ be a bipartite graph satisfying $|E|\leq (1-\gamma)|A||B|$ and $\db(G)\leq k$. Let $G'=(A',B';E')$ be the bipartite graph, where $A'$ is the union of $|B|$ disjoint copies of $A$, $B'$ is the union of $|A|$ disjoint copies of $B$, and there is an edge between copies of vertices if there is an edge between them in $G$. Then $|A'|=|B'|=|A||B|=:n$, and $|E'|=|E||A||B|\leq (1-\gamma)n^2$.
    
    Next, we claim that $\db(G')\leq k$. Say that two vertices of $G'$ are twins if they are copies of the same vertex of $G$. Let $X\subset A'$ and $Y\subset B'$ such that $G'[X,Y]$ is $C_4$-free. If $v,v'\in X\cup Y$ are twins, then $v$ and $v'$ have identical neighbors in $G'[X\cup Y]$, so $C_4$-freeness implies that $v$ and $v'$ have at most one neighbor. Let $X_0\subset X$ be the set of vertices with no twin in $X$, and define $Y_0\subset Y$ analogously. Then $G'[X_0,Y_0]$ is isomorphic to an induced $C_4$-free subgraph of $G$, thus $e(G'[X_0,Y_0])\leq \frac{k}{2}(|X_0|+|Y_0|)$. But then
    \begin{align*}
    e(G'[X,Y])&\leq e(G'[X_0,Y_0])+(|X|-|X_0|)+(|Y|-|Y_0|)\\
    &\leq \frac{k}{2}(|X_0|+|Y_0|)+(|X|-|X_0|)+(|Y|-|Y_0|)\leq \frac{k}{2}(|X|+|Y|),
    \end{align*}
    so $G'[X,Y]$ has average degree at most $k$.

    Therefore, we can apply \Cref{thm:main_EH_balanced} to $G'$ to conclude that there exists an anticomplete pair $(X,Y)$ in $G'$ with $X\subset A'$, $Y\subset B'$ and $|X|,|Y|\geq n^{1-o(1)}$. Let $X_0\subset A$ and $Y_0\subset B$ be the set of vertices which have a copy in $X$ and $Y$, respectively. Then $(X_0,Y_0)$ is anticomplete in $G$, and $|X_0|\geq |X|/|B|$ and $|Y_0|\geq |Y|/|A|$. Thus
    $$|X_0||Y_0|\geq \frac{|X||Y|}{|A||B|}\geq n^{1-o(1)}=(|A||B|)^{1-o(1)},$$
    finishing the proof.
    \end{proof}
    
    The rest of this section is devoted to the proof of \Cref{thm:main_EH_balanced}. In particular, we will focus on \emph{balanced bipartite graphs}, that is, bipartite graphs in which the two vertex classes have equal sizes. First, we prove a technical lemma about finding large subgraphs in a balanced bipartite graph which are somewhat close to being regular. This motivates the following definition.
	
	\begin{definition}
		Let $\beta>0$, and let $G=(A,B;E)$ be a bipartite graph with $|A|=|B|=n$ and maximum degree $D$. Then $G$ is \textbf{$\beta$-semiregular} if for every $X\subset A$ and $Y\subset B$ satisfying $|X|,|Y|\geq n/2$, 
		we have
		$$e(X,Y)\geq \beta Dn.$$
	\end{definition}
	
	\begin{lemma}\label{lemma:semiregular}
		Let $\eps>0$, then the following holds if $n$ is sufficiently large. Let $G=(A,B;E)$ be a bipartite graph with $|A|=|B|=n$, then there exist $X\subset A, Y\subset B$ such that $|X|=|Y|=\Omega_{\eps}(n)$ and $G[X,Y]$ is $n^{-\eps}$-semiregular.  
	\end{lemma}

	\begin{proof}
		Let $\beta=n^{-\eps}$. Let $H=G$, and run the following algorithm, which prunes $H$ in every step. Write $n_H$ for the size of the vertex classes of $H$.
		\medskip
		\begin{algorithmic}[1]
			\State set $H:=G$
			\While{$H$ is not $\beta$-semiregular and $d(H)\geq 1$}
			\If{the maximum degree of $H$ is more than $\beta^{-1/2}d(H)$}
			\State Remove a vertex of maximum degree and an arbitrary vertex from the other class 
			\Else
			\State Let $X,Y$ violate the $\beta$-semiregularity with $|X|=|Y|\geq n_H/2$. Set $H:=H[X,Y]$.
			\EndIf
			\EndWhile
			\State \Return $H$
		\end{algorithmic}
		Let $H_0$ be an instance of $H$ at some step of the algorithm, and write $n_0$ for its size. If $H_0$ is $\beta$-semiregular, or $d(H_0)<1$, then we stop. Otherwise, let $H_1$ be the graph we get in the next step, and let $n_1$ be its size. We consider two cases.
		\begin{description}
			\item[Case 1.] The maximum degree of $H_0$ is more than $\beta^{-1/2}d(H_0)$.
			
			Then line 4 of the algorithm is executed. Let $H_1$ be the graph we get by removing a vertex of maximum degree and  an arbitrary vertex from the other class. Then $n_{1}=n_{0}-1$ and 
			$$e(H_1)\leq d(H_0)n_{0}-\beta^{-1/2}d(H_0)\leq d(H_0) (n_0-\beta^{-1/2}).$$
			Thus,
			$$d(H_1)=e(H_1)n_1^{-1}\leq \frac{d(H_0)(n_0-\beta^{-1/2})}{n_0-1}\leq d(H_0)\left(1-\frac{\beta^{-1/2}}{2n}\right).$$
			
			\item[Case 2.] The maximum degree of $H_0$ is at most $\beta^{-1/2}d(H_0)$.
			
			Then line 6 of the algorithm is executed. Let $X,Y$ violate the $\beta$-semiregularity with $|X|=|Y|\geq n_0/2$. Then $H_1=H_0[X,Y]$, $n_1\geq n_0/2$, and
			$$e(H_1)\leq \beta (\beta^{-1/2}d(H_0)) n_0\leq 2\beta^{1/2} d(H_0) n_1.$$
			Thus, $d(H_1)=e(H_1)n_1^{-1}\leq 2\beta^{1/2}d(H_0).$
		\end{description}
		In both cases, the average degree decreases. In Case 1., the average degree decreases by a factor of $(1-\beta^{-1/2}/(2n))$, so line 4 of the algorithm is executed at most $O(\beta^{1/2}n\log n)=O(n^{1-\eps/2}\log n)$ times, removing a single vertex from both vertex classes. In Case 2., the average degree decreases by a factor of $2\beta^{1/2}$, so line 6 of the algorithm is executed at most $O(1/\eps)$ times. Each such time, the number of vertices decreases by a factor of at most $2$. Thus, the  graph $H$ resulting from the output of our algorithm has at least 
		$$n2^{-O(1/\eps)}-O(n^{1-\eps/2}\log n)\geq \Omega_{\eps}(n).$$
		vertices. If $H$ is $\beta$-semiregular, this $H$ suffices. Otherwise, $d(H)<1$. But then we can select half of the vertices in each class so that there are no edges between them. The resulting empty bipartite graph is also $\beta$-semiregular, finishing our proof.
	\end{proof}
	
	The following definition plays a crucial role in the proof of \Cref{thm:main_EH_balanced}. Recall that if $\cH$ is a hypergraph, then $\cH_h$ denotes  the set of edges of $\cH$ of size exactly $h$.
	
	\begin{definition}
		Let $D,s$ be positive integers, let $L>0$, and let $G=(A,B;E)$ be a bipartite graph with $|A|=|B|=n$.  Then $G$ has the \textbf{$(L,D,s)$-property} if there is a hypergraph $\cH$ on vertex set $A$ such that the edges of $\cH$ have sizes $2,\dots,s$, 
		$$\alpha(\cH[N(y)])\leq s-1\quad\text{for every }y\in B,$$
		and
		$$|\cH_h|\leq LnD^{h-1}\quad\text{for }h=2,\dots,s.$$
	\end{definition}
	
	\noindent
	We highlight a few simple properties of this definition, which will be used in our proofs.
	
	\begin{claim}\label{claim:property}
		Let $G=(A,B;E)$ be a bipartite graph with the $(L,D,s)$-property and  $|A|=|B|=n$.
        \begin{enumerate}\setlength{\itemsep}{0pt}
			\item If $s=1$, then $G$ is the empty graph.
			\item If $D_0<D$, then $G$ has the $(L(D/D_0)^{s-1},D_0,s)$-property.
			\item If $X\subset A$ and $Y\subset B$ such that $|X|=|Y|=n'$,  then $G[X,Y]$ has the $(Ln/n',D,s)$-property.
		\end{enumerate}
	\end{claim}
	
	\begin{proof}
		\begin{enumerate}
			\item If $s=1$, the condition $\alpha(\cH[N(y)])\leq s-1$ implies that the neighborhood of every vertex in $B$ is empty.
			
			\item Writing $L'=L(D/D_0)^{s-1}$, we have $$|\cH_h|\leq LnD^{h-1}=L(D/D_0)^{h-1} n D_0^{h-1}\leq L'n D_0^{h-1}.$$
			
			\item Let $\cH'=\cH[X]$. Then $\cH'$ witnesses the $(Ln/n',D,s)$-property as
			$$|\cH'_h|\leq |\cH_h|\leq LnD^{h-1}\leq L(n/n')\cdot n'D^{h-1}.$$
		\end{enumerate}
	\end{proof}
	
	Our first crucial result shows that bipartite graphs $G$ with $\db(G)\leq k$ have the $$(O_k(1),\Delta(G),k)\text{-property}.$$
	
	\begin{lemma}\label{lemma:k-property}
		Let $k\geq 2$ be an integer, then there is a constant $C>0$ such that the following holds. Let $G=(A,B;E)$ be balanced such that $\db(G)\leq k$ and every vertex in $B$ has degree at most $D$. Then $G$ has the $(C,D,k)$-property.
	\end{lemma}
	
	\begin{proof}
		Let $n=|A|=|B|$. We define $\cH$ by rounding the solution of a fractional version of the problem. Let $\cF$ be the family of all $k$ element subsets of $A$ contained in some neighborhood $N(y)$, $y\in B$.  Moreover, let $\cT$ be the family of all subsets $T\subset A$ such that $2\leq |T|\leq k$. If $\cF$ is empty, simply take $\cH=\emptyset$. Otherwise, consider the linear program
		\begin{align*}
			\text{minimize}\quad &\sum_{T\in \cT}D^{k-|T|}w_T,\\
			\text{subject to}\quad & w_T\geq 0\quad\text{and}\quad\sum_{\substack{T\in \cT\\T\subseteq F}} w_T\geq 1\quad\forall F\in \cF.
		\end{align*}
		Let $(w_T)_{T\in \cT}$ be an optimal solution of this linear program and let 
		$$M=\sum_{T\in \cT}D^{k-|T|}w_T.$$ Then we define $\cH$ to be the hypergraph, whose edges are the sets $T\in \cT$ satisfying $w_T\geq 2^{-k}$. First, we show that $\cH$ satisfies $\alpha(\cH[N(y)])\leq k-1$ for every $y\in B$. Assume that $\alpha(\cH[N(y)])\geq k$ for some $y\in B$, then there exists an independent set $F\subset N(y)$ of size $k$. As $F\in \cF$, we have $\sum_{T\in \cT, T\subseteq F} w_T\geq 1$. As this sum contains $2^k-k-1$ terms, there is some $T\subset F$ such that $w_T\geq 2^{-k}$. But then $T$ is an edge of $\cH$, contradicting that $F$ is an independent set.
		
		Therefore, it remains to show that $|\cH_h|\leq O_k(nD^{h-1})$ for $h=2,\dots,k$. Note that
		$$|\cH_h|\leq 2^k\sum_{T\in \cT, |T|=h} w_T\leq 2^k D^{h-k}M,$$
		thus, it is enough to prove that $M=O_k(nD^{k-1})$. We prove an upper bound on $M$ by considering the dual linear program.
		
		The dual linear program is
		\begin{align*}
			\text{maximize}\quad &\sum_{F\in \cF}\lambda_F,\\
			\text{subject to}\quad & \lambda_F\geq 0\quad\text{and}\quad\sum_{\substack{F\in \cF\\ F\supseteq T}} \lambda_F\leq D^{k-|T|}\quad\forall T\in \cT.
		\end{align*}
		Let $(\lambda_F)_{F\in \cF}$ be an optimal solution of the dual problem, then by the Strong duality theorem, we have
		$$M=\sum_{F\in \cF}\lambda_F.$$
		Observe that $\lambda_F\leq 1$ for every $F\in \cF$, which can be seen by taking $T=F$ in the condition of the program.
		
		Sample every element of $A$ independently with probability $p=1/(2D)$, and let $X$ be the resulting set. Say that $F\in \cF$ is \emph{excellent} if $N(y)\cap X=F$ for some $y\in B$. By fixing any $y\in B$ such that $F\subset N(y)$, we can lower bound the probability of being excellent by
		$$\mathbb{P}(F\text{ is excellent})\geq p^k(1-p)^{|N(y)|-k}\geq p^k (1-p)^D\geq p^k/2.$$
		Also, if $F,F'\in \cF$, we have
		$$\mathbb{P}(F,F'\text{ are excellent})\leq \mathbb{P}(F\cup F'\subset X)=p^{2k-|F\cap F'|}.$$
		Next, sample every member $F\in \cF$ with probability $c_0\lambda_F$, where $c_0>0$ is some sufficiently small constant depending only on $k$, defined later. Let $\cF'$ be the family of sampled sets that are also excellent, then
		$$\mathbb{E}|\cF'|=\sum_{F\in \cF} c_0\lambda_F \mathbb{P}(F\text{ is excellent})\geq \frac{c_0p^k}{2}M.$$
		
		For $h=2,\dots,k$, let $Z_h$ be the number of pairs $\{F,F'\}\in (\cF')^{(2)}$ such that $|F\cap F'|=h$. Then we can bound the expectation of $Z_h$ as 
		\begin{align*}
			\mathbb{E}Z_h&=\sum_{\substack{\{F,F'\}\in \cF^{(2)}\\ |F\cap F'|=h}} c_0^2\lambda_F\lambda_{F'}\mathbb{P}(F,F'\text{ are excellent})\\
			&\leq c_0^2p^{2k-h} \sum_{\substack{\{F,F'\}\in \cF^{(2)}\\ |F\cap F'|=h}} \lambda_F\lambda_{F'}\leq c_0^2p^{2k-h}\sum_{F\in \cF}\lambda_F\sum_{\substack{T\subseteq F\\ |T|=h}}\sum_{\substack{F'\in \cF\\ F'\supseteq T}}\lambda_{F'}\\
			&\leq c_0^2p^{2k-h}\sum_{F\in \cF}\lambda_F\sum_{\substack{T\subseteq F\\ |T|=h}} D^{k-h}\leq c_0^2p^{2k-h}D^{k-h}k^h\sum_{F\in \cF}\lambda_F\\
			&=c_0^2p^{2k-h}D^{k-h}k^hM<c_0^2k^hp^kM
		\end{align*}
		Writing $Z=Z_2+\dots+Z_k$, we thus get 
		$$\mathbb{E}Z\leq 2c_0^2k^kp^kM.$$
		By choosing $c_0$ sufficiently small with respect to $k$, we thus have
		$$\mathbb{E}(|\cF'|-Z)\geq c_1p^kM$$
		for some $c_1>0$ depending only on $k$. For any pair $F,F'\in \cF'$ satisfying $|F\cap F'|\geq 2$, remove one of the sets from $\cF'$, let $\cF''$ be the resulting family. For every $F\in \cF''$, choose some vertex $y\in B$ such that $X\cap N(y)=F$, and let $Y$ be the set of selected vertices. Then $|Y|=|\cF''|$, and thus $$\mathbb{E}|Y|=\mathbb{E}|\cF''|\geq c_1p^kM.$$ On the other hand $G'=G[X,Y]$ is a $C_4$-free induced subgraph of $G$. Indeed, if $y,y'\in Y$ are distinct, then $N_{G'}(y)=F$ and $N_{G'}(y')=F'$ for some distinct $F,F'\in \cF''$, so the condition $|F\cap F'|<2$ ensures that $y$ and $y'$ are not the vertices of a four-cycle. Hence, $\db(G)\leq k$ implies
		$$e(G[X,Y])\leq \frac{k}{2}(|X|+|Y|).$$
		On the other hand, $e(G[X,Y])=k|Y|$, as every vertex of $Y$ has exactly $k$ neighbors. Therefore, the previous inequality is equivalent to $|Y|\leq |X|$. Taking expectations, we  get
		$$c_1p^kM\leq \mathbb{E}|Y|\leq \mathbb{E}|X|\leq pn,$$
		which implies
		$$M\leq 2^{k-1}c_1^{-1}nD^{k-1}.$$
		This finishes the proof.
	\end{proof}
	
	Next, we prove that if a bipartite graph $G$ has the $(L,D,s+1)$-property, and is sufficiently semiregular, then it contains a large induced subgraph with $(L',D,s)$-property for some well controlled $L'$. Given such a result, we prove \Cref{thm:main_EH_balanced} by iterating this step: starting with a bipartite graph having $(O_k(1),D,k)$-property, which is guaranteed by the previous lemma, we find a large induced subgraph with the $(L',D,1)$-property for some $L'$. But then the vertex classes of this bipartite  subgraph form an anticomplete pair, finishing the proof. Executing this approach formally has a number of difficulties, so the precise argument is postponed.

	\begin{lemma}\label{lemma:s-property}
		Let $s\in \mathbb{Z}^+$, then there exists $c_s>0$ such that the following holds. Let $n,D$ be positive integers and $L\geq 1$ such that $128LD\leq  n$. Let $G=(A,B;E)$ be a $\beta$-semiregular bipartite graph with $|A|=|B|=n$ and maximum degree $D$. Assume that $G$ has the $(L,D,s+1)$-property.
		Then there are sets $X\subset A$, $Y\subset B$ such that $|X|=|Y|\geq \Omega(\beta n/L)$, and $G[X,Y]$ has the $(c_s L^3/\beta^2,D,s)$-property.
	\end{lemma}
	
	\begin{proof}
		Let $\cH$ be a hypergraph on $A$ witnessing that $G$ has the $(L,D,s+1)$-property. First, we remove all vertices from $A$, whose degree in $\cH_2$ is more than
		$$\Delta:=n/32.$$
		As $|\cH_2|\leq LnD$, we removed at most 
		$$\frac{2LnD}{\Delta}=64LD\leq \frac{n}{2}$$
		vertices. Let $A'$ be the resulting set. For $x\in A'$, let 
		$$Q_x=\{v\in A': \{x,v\}\in \cH_2\}$$
		be the neighbors of $x$ in $\cH_2$, then $|Q_x|\leq \Delta$. Set 
		$$T=\sum_{x\in A'}|Q_x|\leq 2|\cH_2|\leq 2LnD.$$
		Also, we note that as $G$ is $\beta$-semiregular, we have $e(A',B)\geq \beta Dn$.
		
		Sample each vertex of $A'$ with probability $p=1/(64D)$, and let $S$ be the set of sampled vertices. Moreover, further sample each element of $S$ with probability 
		$q=1/(2L),$ let $S^*$ be the resulting set. Define $$X=A'\setminus A_0\quad\text{with }A_0=S\cup \bigcup_{x\in S^*}Q_x,$$
		and
		$$Y=\{y\in B: |N(y)\cap S|=1\text{ and } N(y)\cap S\subset S^*\}.$$
		First, we show that $X$ and $Y$ are simultaneously large with some non-negligible probability. Write $$\lambda=\frac{\beta}{2048L}.$$
		
		\begin{claim}
			We have
			$$\mathbb{P}(|X|\geq n/4\text{ and }|Y|\geq \lambda n)\geq \lambda.$$
		\end{claim}
		\begin{proof}
			For $y\in B$, let $d_y=|A'\cap N(y)|\leq D$, then
			$$\mathbb{P}(y\in Y)=q\cdot d_yp(1-p)^{d_y-1}\geq \frac{qd_y}{128D}.$$
			Fix some $y\in B$ and $z\in N(y)\cap A'$, and consider the event 
			$$\mathcal{E}_{y,z}=\{N(y)\cap S=\{z\}\text{ and }z\in S^*\}.$$ 
			Then, recalling that $A_0=S\cup \bigcup_{x\in S^*}Q_x$, we have
			\begin{align*}
				\mathbb{E}(|A_0| \mid \mathcal{E}_{y,z})&\leq \mathbb{E}(|S| \mid \mathcal{E}_{y,z})+\mathbb{E}\left(\left|\bigcup_{x\in S^*}Q_x\right| \mid \mathcal{E}_{y,z}\right)\\
				&\leq 1+\mathbb{E}(|S|)+|Q_z|+\mathbb{E}\left(\left|\bigcup_{x\in S^*}Q_x\right|\right)\\
				&\leq  1+pn+\Delta+pqT\leq 1+pn+n/32+n/64<n/8.
			\end{align*}
			Therefore, by Markov's inequality, $\mathbb{P}(|A_0|>n/4  \mid \mathcal{E}_{y,z})\leq 1/2$. Recalling that $X=A'\setminus A_0$ and $|A'|\geq n/2$, this implies that
			$$\mathbb{P}(|X|< n/4\mid  \mathcal{E}_{y,z})\leq \mathbb{P}(|A_0|>n/4  \mid \mathcal{E}_{y,z})\leq 1/2.$$
			Note that $y\in Y$ if and only if $\mathcal{E}_{y,z}$ holds for some $z\in N(y)$. Therefore, the previous implies
			$$\mathbb{P}(|X|\geq n/4\mid y\in Y)\geq 1/2,$$
			from which
			$$\mathbb{P}(|X|\geq n/4\text{ and }y\in Y)\geq \frac{1}{2}\mathbb{P}(y\in Y)\geq \frac{qd_y}{256D}.$$
			Therefore,
			$$\mathbb{E}(|Y|\mathds{1}_{|X|\geq n/4})=\sum_{y\in B}\mathbb{P}(|X|\geq n/4\text{ and }y\in Y)\geq \frac{q}{256D}\sum_{y\in B}d_y=\frac{qe(A',B)}{256D}\geq \frac{\beta n}{512L}=4\lambda n.$$
			As $|Y|\mathds{1}_{|X|\geq n/4}\leq n$, \Cref{lemma:simple_prob} applied with $c=1/2$ implies that 
			$$\mathbb{P}(|Y|\mathds{1}_{|X|\geq n/4}\geq 2\lambda n)\geq \frac{4\lambda n}{2n}\geq \frac{\beta}{1024L}.$$
			This finishes the proof of the claim.
		\end{proof}

		Define the hypergraph $\cH'$ on $X$ as follows. For $h=2,\dots,s$, let
		$$\cH'_h=\cH_h[X]\cup \{F\setminus \{x\}: F\in \cH_{h+1}, x\in F\cap S^*, F\setminus\{x\}\subset X\}.$$ 
		That is, we take all edges of $\cH[X]$ of size $h$, and in addition the truncations of those edges of size $h+1$ in $\cH$, that intersect $A_0$ in a single element of $S^*$. For every $F\in \cH_{h+1}$, the probability that some $x\in F$ is contained in $S^*$ is at most $(h+1)pq\leq (s+1)/D$. Thus, we have
		$$\mathbb{E}|\cH'_h|\leq |\cH_h|+\frac{s+1}{D}|\cH_{h+1}|\leq (s+2)LnD^{h-1}.$$
		By Markov's inequality, this gives
		$$\mathbb{P}\left(|\cH'_h|\geq \frac{s(s+2)}{\lambda}LnD^{h-1}\right)\leq \frac{\lambda}{s}.$$
		Therefore, with positive probability, the events
		$$|X|\geq \lambda n,\quad |Y|\geq \lambda n,\quad |\cH'_h|\leq \frac{s(s+2)}{\lambda}LnD^{h-1}\quad (2\leq h\leq s)$$
		hold simultaneously. Thus, we can fix an outcome of $S$ and $S^*$ such that all these events hold.
		
		Finally, we prove that for every $y\in Y$, $\cH'[N(y)\cap X]$ has no independent set of size $s$. Recall that if $y\in Y$, then there is a unique $\{z\}= S^*\cap N(y)$ for some vertex $z$. Let $J\subset N(y)\cap X$ of size $s$, and assume to the contrary that $J$ is an independent set in $\cH'$. As $\cH'$ contains all edges of $\cH[X]$ of size at most $s$, $J$ must be an independent set in $\cH$ as well. But as $J\cup \{z\}$ is not an independent set in $\cH[N(y)]$, there exists some edge $F$ in $\cH$ such that $F\subset J\cup\{z\}$. But then we must have $z\in F$. If $F$ has size 2, then $F=\{v,z\}\in \cH_2$ for some $v\in Q_z\subset A_0$. But this is impossible as $X$ is disjoint from $Q_z$, contradicting $v\in X$. On the other hand, if $|F|\geq 3$, then $F\setminus\{z\}$ is an edge of $\cH'$, and $F\setminus\{z\}\subset J$, contradicting that $J$ is an independent set. 
		
		Let $n'=\min\{|X|,|Y|\}$. We may remove elements arbitrarily from the larger set to make them of equal size, without violating any of the necessary conditions. Thus, we have $n'=|X|=|Y|\geq \lambda n$, and for $h=2,\dots,s$,
		$$|\cH_h'|\leq \frac{s(s+2)}{\lambda^2}Ln'D^{h-1}=O_s(L^3/\beta^2) n'D^{h-1}.$$
		Thus, $G[X,Y]$ has the $(O_s(L^3/\beta^2),D,s)$-property. This finishes the proof.
	\end{proof}
	
	Using the previous lemmas, we prove the special subcase of \Cref{thm:main_EH_balanced} in which we have a strong upper bound on the maximum degree. The proof of the general case then will be a simple reduction to this case by standard VC-dimension arguments.
	
	\begin{lemma}\label{lemma:max_deg}
		Let $k\geq 2$ be an integer, $\tau>0$, then the following holds if $n$ is sufficiently large. Let $G=(A,B;E)$ be a bipartite graph such that $|A|=|B|=n$, $\db(G)\leq k$, and the maximum degree of $G$ is at most $n^{1-\tau}$. Then $G$ contains an anticomplete pair of size at least $n^{1-\tau}$.
	\end{lemma}
	
	\begin{proof}
		Let $C$ be the constant given by \Cref{lemma:k-property}, depending only on $k$, and let $c_s$ be the constant given by \Cref{lemma:s-property}, depending only on $s$. Choose $\eps>0$ sufficiently small with respect to $\tau$ and $k$, specified later. Furthermore, let $\rho=\rho(\eps)>0$ be the constant given by \Cref{lemma:semiregular}. We assume that $n$ is sufficiently large with respect to $\eps,k,\tau$, which will be implicitly used in our calculations.
		
		Let $\beta=n^{-\eps}$, and define two sequences of positive real numbers $L_1,\dots,L_k$ and $x_2,\dots,x_k$ as follows. Let $L_k=C$, and for $s=k,k-1,\dots,2$, let 
		$$x_s=(\beta^{-k}L_k\dots L_s)^{2/\tau}\quad\text{and}\quad L_{s-1}=c_{s-1}L_s^3x_s^{s-2}\beta^{-2}\rho^{-1}.$$
		Clearly, we have $x_s,L_s\leq n^{O_{k,\tau}(\eps)}$. 
		
		By \Cref{lemma:semiregular}, we can choose a
		$\beta$-semiregular induced subgraph $G_k$ of $G$
		with at least $\rho n$ vertices in both vertex classes. We run the following algorithm, whose input is the graph $G_k$.
		
		\begin{enumerate}
			\item \textbf{Initialize.}
			If $G_k$ is edgeless, stop and output $G_k$.
			Otherwise, $G_k$ satisfies the conditions below.
			
			\item \textbf{Iterate.}
			For $s=k,k-1,\ldots,2$, assume that $G_s$ has
			$n_s$ vertices in both vertex classes, maximum degree $D_s$, is
			$\beta$-semiregular, and has the
			$(L_s,D_s,s)$-property. Perform the following steps.
			
			\begin{enumerate}
				\item If $128L_sD_s>n_s$, then stop and output FAILED.
				\item 
				By \Cref{lemma:s-property}, we can choose $H=G_s[X,Y]$
				such that $|X|=|Y|=\Omega(\beta n_s/L_s)$ and $H$ has the
				$(M_s,D_s,s-1)$-property, where
				\[
				M_s=\frac{c_{s-1}L_s^3}{\beta^2}.
				\]
				If $s=2$, stop and output $H$.
				
				\item 
				By \Cref{lemma:semiregular}, choose a
				$\beta$-semiregular induced subgraph $G_{s-1}$
				of $H$ with vertex classes of size $n_{s-1}\ge \rho|X|$.
				By item~3 of \Cref{claim:property}, it has the
				$(M_s/\rho,D_s,s-1)$-property.
				If $G_{s-1}$ is edgeless, stop and output it.
				Set $D_{s-1}=\Delta(G_{s-1})$.
				
				\item 
				If $x_sD_{s-1}<D_s$, set $G_k\gets G_{s-1}$
				and \emph{restart} at step~1.
				
				Otherwise, $D_s/D_{s-1}\le x_s$, so item~2 of
				\Cref{claim:property} gives the $(x_s^{s-2}M_s/\rho,D_{s-1},s-1)$-property.
				By the definition of $L_{s-1}$, this is equal to the $(L_{s-1},D_{s-1},s-1)$-property. The required conditions are satisfied
				for the next iteration.
			\end{enumerate}
		\end{enumerate}

		We now analyze this algorithm. We first look at how the sizes and maximum degrees change within one \emph{cycle} of the algorithm, which refers to the steps executed before the algorithm is restarted.
		Observe that
		$$\frac{n_{k}}{n_{s-1}}\leq O_{k}\left((\rho\beta)^{-(k-s+1)}L_k\dots L_s\right)< x_s^{\tau/2}.$$ 
		Therefore, if we restarted the cycle at some $s$, then we have 
		$$\frac{n_{k}}{n_{s-1}}\leq \left(\frac{D_k}{D_{s-1}}\right)^{\tau/2}.$$
		In particular, this implies the following. Let $G_s^{(\ell)}$ denote the graph $G_s$ in the $\ell$-th cycle, and write $n_s^{(\ell)}=n_s$, $D_s^{(\ell)}=D_s$. Then
		$$\frac{n_k^{(\ell)}}{n_k^{(\ell+1)}}\leq \left(\frac{D_k^{(\ell)}}{D_k^{(\ell+1)}}\right)^{\tau/2},$$
		which implies by taking the telescopic product that
		$$\frac{n_{k}^{(1)}}{n_k^{(\ell)}}\leq \left(\frac{D_k^{(1)}}{D_k^{(\ell)}}\right)^{\tau/2}\leq n^{\tau/2}.$$
		Hence, as $n_k^{(1)}\geq \rho n$, each cycle starts with a balanced bipartite graph whose vertex classes have size at least $n^{1-2\tau/3}$. Moreover, as $n_k^{(\ell)}/n_s^{(\ell)}\leq n^{O(\eps)}$, by choosing $\eps$ sufficiently small, we also have $n_s^{(\ell)}\geq n^{1-3\tau/4}$ for every valid $s$ and $\ell$. But then as $D_s^{(\ell)}\leq n^{1-\tau}$ and $L_s=n^{O_{k,\tau}(\eps)}$, we always have $128L_s^{(\ell)}D_s^{(\ell)}\leq n_s^{(\ell)}$ satisfied (assuming $\eps$ is sufficiently small with respect to $k$ and $\tau$), so the algorithm cannot output FAILED. But then the algorithm outputs an anticomplete pair of size at least $n^{1-\tau}$.
	\end{proof}
	
	Finally, we prove \Cref{thm:main_EH_balanced}. It only remains to show that every balanced degree-bounded bipartite graph with density at most $1-\gamma$ contains a large induced subgraph with small maximum degree. This is a standard consequence of the \emph{Ultra-strong regularity lemma} for graphs of bounded VC-dimension. Given a bipartite graph $G=(A,B;E)$ with $|A|=|B|=n$, an \emph{$r$-equipartition} of $G$ is a partition of $A$ into $r$ parts $A_1,\dots,A_r$, and a partition of $B$ into $r$ parts $B_1,\dots,B_r$, such that $|A_i|,|B_i|\in \{\lfloor n/r\rfloor,\lceil n/r\rceil\}$. Such an equipartition is $\delta$-homogeneous if there are at most $\delta r^2$ pairs $(i,j)\in [r]^2$ such that the edge density of $G[A_i,B_j]$ is in the interval $(\delta,1-\delta)$. The following is a simple corollary of Lemma 5.6 in \cite{AFN}.
	
	\begin{lemma}\label{lemma:homogeneous}
		Let $\delta\in (0,1)$, and let $G=(A,B;E)$ be a balanced bipartite graph which does not contain a fixed bipartite graph $H$ as an induced subgraph. If $H$ has $q$ vertices, then $G$ has a $\delta$-homogeneous $r$-equipartition for some $r<(q/\delta)^{O(q)}$. 
	\end{lemma}
	
	\begin{proof}[Proof of \Cref{thm:main_EH_balanced}]
		Let $\tau>0$ be sufficiently small with respect to $\eps$ and $\gamma$, specified later. It is well known that there exists a bipartite $C_4$-free graph $H$ on $q=O(k^2)$ vertices with average degree more than $k$: one may take a point-line incidence graph over the finite projective plane of order $p=\Theta(k)$. Then $G$ does not contain $H$ as an induced subgraph. Thus, applying \Cref{lemma:homogeneous} with $\delta=n^{-2\tau}$, we get a $\delta$-homogeneous $r$-equipartition of $G$ for some $r<(q/\delta)^{O(q)}<n^{\alpha}/3$, where $\alpha=Cq\tau$ for some $C>0$ sufficiently large. Let $A_1,\dots,A_r$, and $B_1,\dots,B_r$ be the parts in this equipartition.
		
		Observe that there exists a pair $(i_0,j_0)\in [r]^2$ such that the density of $G[A_{i_0},B_{j_0}]$ is at most $\delta$. Otherwise, there are at least $(1-\delta)r^2$ pairs $(i,j)$ such that the density of $G[A_i,B_j]$ is at least $1-\delta$, which implies that $G$ has at least $(1-\delta)r^2 (1-\delta)\lfloor n/r\rfloor^2>(1-\gamma) n^2$ edges, contradiction. Remove all vertices of degree more than $4\delta \lceil n/r\rceil$ in $G[A_{i_0},B_{j_0}]$, then we get two sets $X\subset A_{i_0}$ and $Y\subset B_{j_0}$  such that $|X|,|Y|\geq n/(2r)$ and $G[X,Y]$ has maximum degree at most $8\delta n/r$. By removing further vertices, we can ensure that $|X|=|Y|\geq n/(2r)$, and set $G_0=G[X,Y]$. Then $G_0$ is a balanced induced subgraph of $G$ of size $n_0\geq n/2r\geq n^{1-\alpha}$ and maximum degree at most $8\delta n/r\leq n_0^{1-\tau}$.
		
		Applying \Cref{lemma:max_deg} for $G_0$ instead of $G$, we find an anticomplete pair of size at least $n_0^{1-\tau}\geq n^{(1-\alpha)(1-\tau)}\geq n^{1-\eps}$, finishing the proof.
	\end{proof}
	
	\section{Coloring degree-bounded graphs}\label{sect:chromatic}
	
	In this section, we prove \Cref{thm:main_color}, which we restate here in the following equivalent form.
    
    \begin{theorem}\label{thm:main_color_precise}
		Let $k,r\geq 2$ be positive integers and $\beta>0$, then the following holds for every sufficiently large $n$. Let $G$ be an $n$-vertex $K_r$-free graph such that $\db(G)\leq k$. Then $\chi(G)\leq n^{\beta}$. 
	\end{theorem}
	
	The proof follows a very similar strategy as the proof of \Cref{thm:main_EH_balanced}. First, we redefine semiregularity so it applies to graphs instead of bipartite graphs.
	
	\begin{definition}
		Let $\beta,c>0$, and let $G$ be an $n$-vertex graph of maximum degree $D$. Then $G$ is \textbf{$(\beta,c)$-semiregular} if for every $X\subset V(G)$ satisfying $|X|\geq cn$, 
		we have
		$$e(X)\geq \beta Dn.$$
	\end{definition}
	
	We will apply $(\beta,c)$-semiregularity only with the values $c=1/4$ or $1/2$, and we remark that  if $G$ is $(\beta,1/4)$-semiregular, then every induced subgraph of $G$ on at least $|V(G)|/2$ vertices is $(\beta,1/2)$-semiregular. Moreover, we prove the following analogue of \Cref{lemma:semiregular}. 
	
	\begin{lemma}\label{lemma:semiregular_nonbip}
		Let $\eps,c>0$, then there exists $\rho>0$ such that the following holds if $n$ is sufficiently large. Let $G$ be an $n$-vertex graph, then there exists $X\subset V(G)$ such that $|X|\geq \rho n$ and $G[X]$ is  $(n^{-\eps},c)$-semiregular.
	\end{lemma}

	\begin{proof}
		Let $\beta=n^{-\eps}$. Let $H=G$, and run the following algorithm, which prunes $H$ in every step. Write $n_H$ for the size of $H$.
		\medskip
		\begin{algorithmic}[1]
			\State set $H:=G$
			\While{$H$ is not $(\beta,c)$-semiregular and $d(H)\geq 1$}
			\If{the maximum degree of $H$ is more than $\beta^{-1/2}d(H)$}
			\State Remove a vertex of maximum degree
			\Else
			\State Let $X\subset V(H)$ violate the $(\beta,c)$-semiregularity with $|X|\geq cn_H$. Set $H:=H[X]$.
			\EndIf
			\EndWhile
			\State \Return $H$
		\end{algorithmic}
		Let $H_0$ be an instance of $H$ at some step of the algorithm, and write $n_0$ for its size. If $H_0$ is $(\beta,c)$-semiregular, or $d(H_0)<1$, then we stop. Otherwise, let $H_1$ be the graph we get in the next step, and let $n_1=|V(H_1)|$. We consider two cases.
		\begin{description}
			\item[Case 1.] The maximum degree of $H_0$ is more than $\beta^{-1/2}d(H_0)$.
			
			Then line 4 of the algorithm is executed. Let $H_1$ be the graph we get by removing a vertex of maximum degree. Then $n_{1}=n_{0}-1$ and 
			$$2e(H_1)\leq d(H_0)n_{0}-\beta^{-1/2}d(H_0)\leq d(H_0) (n_0-\beta^{-1/2}).$$
			Thus,
			$$d(H_1)=2e(H_1)n_1^{-1}\leq \frac{d(H_0)(n_0-\beta^{-1/2})}{n_0-1}\leq d(H_0)\left(1-\frac{\beta^{-1/2}}{2n}\right).$$
			
			\item[Case 2.] The maximum degree of $H_0$ is at most $\beta^{-1/2}d(H_0)$.
			
			Then line 6 of the algorithm is executed. Let $X$ violate the $(\beta,c)$-semiregularity with $|X|\geq cn_0$. Then $H_1=H_0[X]$, $n_1\geq cn_0$, and
			$$2e(H_1)\leq 2\beta \Delta(H_0)n_0\leq 2\beta^{1/2}d(H_0) n_0\leq \frac{2\beta^{1/2}}{c} d(H_0) n_1.$$
			Thus, $d(H_1)=2e(H_1)n_1^{-1}\leq 2\beta^{1/2}d(H_0)/c.$
		\end{description}
		In both cases, the average degree decreases. In Case 1., the average degree decreases by a factor of $(1-\beta^{-1/2}/(2n)$, so line 4 of the algorithm is executed at most $O(\beta^{1/2}n\log n)=O(n^{1-\eps/2}\log n)$ times, removing one vertex. In Case 2., the average degree decreases by a factor of $2\beta^{1/2}/c$, so line 6 of the algorithm is executed at most $O_c(1/\eps)$ times. Each such time, at least $c$-fraction of the vertices is retained. Thus, the  graph $H$ resulting from the output of our algorithm has at least 
		$$nc^{O_c(1/\eps)}-O(n^{1-\eps/2}\log n)\geq \Omega_{\eps,c}(n).$$
		vertices. If $H$ is $(\beta,c)$-semiregular, this $H$ suffices. Otherwise, $d(H)<1$. But then we can select half of the vertices so that the resulting graph is empty, finishing our proof.
	\end{proof}

    Next, we redefine the $(L,D,s)$-property. We define the analogue, called stable-$(L,D,s)$-property, which has the additional requirement that the edges of the mixed hypergraph $\cH$ are independent sets in the graph.
	
	\begin{definition}
		Let $D,s$ be positive integers, let $L>0$, and let $G$ be an $n$-vertex graph.  Then $G$ has the \textbf{stable-$(L,D,s)$-property} if there is a hypergraph $\cH$ on vertex set $V(G)$ such that the edges of $\cH$ have sizes $2,\dots,s$, every edge of $\cH$ is an \emph{independent set} in $G$, 
		$$\alpha(\cH[N(y)])\leq s-1\quad\text{for every }y\in V(G),$$
		and
		$$|\cH_h|\leq LnD^{h-1}\quad\text{for }h=2,\dots,s.$$
	\end{definition}
	
	Akin to \Cref{claim:property}, we collect some simple properties of this definition, which will be used later.
	
	\begin{claim}\label{claim:stable_property}
		Let $G$ be an $n$-vertex graph with the stable-$(L,D,s)$-property. 
		\begin{enumerate}\setlength{\itemsep}{0pt}
			\item If $s=1$, then $G$ is the empty graph.
			\item If $D_0<D$, then $G$ has the stable-$(L(D/D_0)^{s-1},D_0,s)$-property.
			\item If $G'$ is a subgraph of $G$ on $n'$ vertices,  then $G'$ has  the stable-$(Ln/n',D,s)$-property.
		\end{enumerate}
	\end{claim}
	
	\begin{proof}
		\begin{enumerate}
			\item If $s=1$, the condition $\alpha(\cH[N(y)])\leq s-1$ implies that the neighborhood of every vertex is empty.
			
			\item Writing $L'=L(D/D_0)^{s-1}$, we have $$|\cH_h|\leq LnD^{h-1}=L(D/D_0)^{h-1} n D_0^{h-1}\leq L'n D_0^{h-1}.$$
			
			\item Let $\cH'=\cH[V(G')]$. Then $\cH'$ witnesses the stable-$(Ln/n',D,s)$-property as
			$$|\cH'_h|\leq |\cH_h|\leq LnD^{h-1}\leq L(n/n')\cdot n'D^{h-1}.$$
		\end{enumerate}
	\end{proof}

    Our first main auxiliary lemma is an analogue of \Cref{lemma:k-property}. We show that a $K_r$-free degree-bounded graph contains a large induced subgraph with the stable-$(L,\Delta(G),O(1))$-property, where $L$ depends on the chromatic number of neighborhoods. Thus, the main difference compared to \Cref{lemma:k-property} is that we are no longer able to choose $L$ to be a constant. However, the neighborhoods are $K_{r-1}$-free, so assuming we proceed by induction on $r$ in the final proof of \Cref{thm:main_color_precise}, we will have $L=n^{o(1)}$ in our applications.
	
	\begin{lemma}\label{lemma:stable}
		Let $k,r\geq 2$ be integers, then there exist $C>0$ and $t\in \mathbb{Z}^+$ such that the following holds. Let $G$ be a $K_r$-free $n$-vertex graph such that $\db(G)\leq k$ and $\Delta(G)\leq D$. Assume that $G[N(x)]$ is $K$-colorable for every $x\in V(G)$. Then there exists $A\subset V(G)$ of size at least $3n/4$ such that $G[A]$ has the stable-$(CK^k,D,t)$-property.
	\end{lemma}
	
	\begin{proof}
		First, we define some parameters depending only on $k$ and $r$. Let $t$ be the smallest number such that every graph on $t$ vertices contains either a clique of size $r$, or an independent set of size $k$. Let $\mathcal{T}$ be the family of all sets $T\subset V(G)$ such that $2\leq |T|\leq k$ and $T$ is an independent set in $G$. Moreover, let $\cF$  be the family of all pairs $(y,F)$ where $y\in V(G)$ and $F\subset N(y)$ is an independent set of size $k$. Let $\Delta=C_1K^kD^{k-1}$ for some $C_1$ depending only on $k$, specified later. Consider the linear program
		\begin{align*}
			\text{minimize}\quad &\sum_{T\in \cT}D^{k-|T|}w_T+\Delta\sum_{y\in V(G)} u_y,\\
			\text{subject to}\quad & w_T\geq 0, u_y\geq 0,\text{ and}\\
			\quad & u_y+\sum_{\substack{T\in \cT\\T\subseteq F}} w_T\geq 1\quad\forall (y,F)\in \cF.
		\end{align*}
		Let $M$ be the optimal value of this linear program and let $(w_T)_{T\in \cT}, (u_y)_{y\in V(G)}$ satisfy the conditions and achieve
		$$M=\sum_{T\in \cT}D^{k-|T|}w_T+\Delta\sum_{y\in V(G)} u_y.$$
		
		Our goal is to prove that $M\leq \Delta n/8$. First, we show that this indeed implies the existence of the desired set $A$.
		
		\begin{claim}
			If $M\leq \Delta n/8$, then there exists a set $A\subset V(G)$ such that $|A|\geq 3n/4$ and $G[A]$ has the stable-$(O_k(K^k),D,t)$-property.
		\end{claim}
		\begin{proof}
			Let
			$$A=\{y\in V(G): u_y\leq 1/2\},$$
			and define the hypergraph $\cH$ on $A$ such that
			$$E(\cH)=\{T\in \cT: T\subset A\text{ and }w_T\geq  2^{-k-1}\}.$$
			We show that $\cH$ witnesses the stable-$(O_k(K^k),D,t)$-property of $G[A]$. Firstly, we have
			$$n-|A|=|\{y\in V(G):u_y>1/2\}|\leq 2\sum_{y\in V(G)}u_y\leq \frac{2M}{\Delta}\leq \frac{n}{4},$$
			proving that $|A|\geq 3n/4$. Secondly, for $h=2,\dots,k$, we have
			$$|\cH_h|\leq 2^{k+1}\sum_{T\in \cT, |T|=h}w_T\leq 2^{k+1}\frac{M}{D^{k-h}}\leq \frac{2^{k}\Delta n}{D^{k-h}}=2^kC_1D^{h-1}K^kn=O_k(D^{h-1}K^k|A|),$$
			and  $|\cH_h|=0$ for $h=k+1,\dots,t$. Thirdly, let $y\in A$ and let $J\subset N(y)\cap A$ be some set of size $t$. Then by the definition of $t$ and using that $G$ is $K_r$-free, there exists $F\subset J$ of size $k$ that is independent in $G$. But as 
			$$u_y+\sum_{T\in\cT, T\subseteq F}w_T\geq 1$$
			and $u_y\leq 1/2$, we have
			$$\sum_{T\in\cT, T\subseteq F}w_T\geq 1/2.$$
			Thus, as there are $2^k$ subsets of $F$, there exists some $T\in \cT$, $T\subset F$, such that $w_T\geq 2^{-k-1}$. But then $T$ is an edge of $\cH$ and $T\subset J$, so $J$ is not an independent set in $\cH[N(y)\cap A]$. This concludes the proof that $G[A]$ has the stable-$(O_k(K^k),D,t)$-property.
		\end{proof}
		
		The rest of the proof is devoted to establishing the inequality $M\leq \Delta n/8$. In order to upper bound $M$, we consider the dual linear program. For every $(y,F)\in \cF$, let $\lambda_{y,F}$ be a variable. Then the dual linear program is 
		\begin{align*}
			\text{maximize}\quad &\sum_{(y,F)\in \cF}\lambda_{y,F},\\
			\text{subject to}\quad & \lambda_{y,F}\geq 0,\\
			\quad & \sum_{(y,F)\in \cF:T\subseteq F} \lambda_{y,F}\leq D^{k-|T|}\quad\forall T\in \cT\text{, and}\\
			\quad & \sum_{F:(y,F)\in \cF}\lambda_{y,F}\leq \Delta\quad \forall y\in V(G)
		\end{align*}
		By the Strong duality theorem, the solution of this program is also $M$, so fix some $(\lambda_{y,F})_{(y,F)\in \cF}$ satisfying the conditions and
		$$M=\sum_{(y,F)\in \cF}\lambda_{y,F}.$$
		We note that by taking $T=F$, the condition $\sum_{(y,F)\in \cF:T\subseteq F}\lambda_{y,F}\leq 1$ ensures that $\lambda_{y,F}\leq 1$ for every $(y,F)\in \cF$.
		
		Let $c_0=1/(1000k)$ and $p=c_0/D$. Sample the vertices of $V(G)$ independently with probability $p$, and let $X$ be the resulting set. Let $S$ be the set of isolated vertices in $G[X]$.  Say that $(y,F)$ is \emph{excellent} if $F\subset X$ and $X$ is disjoint from 
		$$(N(y)\setminus F)\cup \bigcup_{x\in F}N(x).$$
		The displayed set has size at most $(k+1)D$, so
		$$\mathbb{P}((y,F)\text{ is excellent})\geq p^k(1-p)^{(k+1)D}\geq \frac{1}{2}p^k.$$
		Observe that if $(y,F)$ is excellent, then every vertex of $F$ is isolated in $G[X]$, so $F\subset S$. Moreover, for every $y$, there is at most one set $F$ such that $(y,F)$ is excellent, as $F=N(y)\cap X$. We write $F_y$ for this unique $F$. Furthermore, if $(y,F),(y',F')\in \cF$, then
		$$\mathbb{P}((y,F),(y',F')\text{ are excellent})\leq \mathbb{P}(F,F'\subset X)=p^{2k-|F\cap F'|}.$$
		
		Next, let $c_1=k^{-k-1}/4$, and sample every $(y,F)\in \cF$ independently with probability $c_1\lambda_{y,F}$. Let $\cF'$ be the set of sampled pairs that are also excellent. Then
		$$\mathbb{E}|\cF'|=\sum_{(y,F)\in \cF}c_1\lambda_{y,F}\mathbb{P}((y,F)\text{ is excellent})\geq \frac{1}{2}c_1p^kM.$$
		For $h=2,\dots,k$, let  $Z_h$ be the number of ordered pairs of distinct members $(y,F),(y',F')\in \cF'$ such that $|F\cap F'|=h$. Then
		\begin{align*}
			\mathbb{E}Z_h&=\sum_{\substack{(y,F)\neq (y',F')\in \cF\\ |F\cap F'|=h}}\mathbb{P}((y,F),(y',F')\text{ are excellent})c_1^2\lambda_{y,F}\lambda_{y',F'}\\
			&\leq c_1^2p^{2k-h}\sum_{\substack{(y,F),(y',F')\in \cF\\ |F\cap F'|=h}}\lambda_{y,F}\lambda_{y',F'}\\
			&\leq c_1^2p^{2k-h}\sum_{(y',F')\in \cF}\lambda_{y',F'}\sum_{\substack{T\subseteq F'\\|T|=h}}\sum_{\substack{(y,F)\in \cF\\ T\subseteq F}}\lambda_{y,F}\\
			&\leq  c_1^2p^{2k-h}\sum_{(y',F')\in \cF}\lambda_{y',F'}\sum_{\substack{T\subseteq F'\\|T|=h}}D^{k-h}\\
			&\leq   c_1^2p^{2k-h}\binom{k}{h}D^{k-h}\sum_{(y',F')\in \cF}\lambda_{y',F'}\\
			&=c_1^2p^{2k-h}\binom{k}{h}D^{k-h}M\leq k^h c_1^2p^kM.
		\end{align*}
		Let $Z=Z_2+\dots+Z_k$, then 
		$$\mathbb{E}Z\leq k^{k+1}c_1^2 p^kM\leq \frac{1}{2}\mathbb{E}|\cF'|.$$ 
		For each pair $(y,F),(y',F')\in \cF'$ satisfying $|F\cap F'|\geq 2$, remove one of them from $\cF'$, let $\cF''$ be the resulting family. Moreover, let $Y_0$ be the set of vertices $y$ such that $(y,F_y)\in \cF''$. Then $|Y_0|=|\cF''|$, and thus $$\mathbb{E}|Y_0|\geq \frac{1}{4}c_1p^k M.$$
		
		For every $x\in V(G)$, let $\chi_x:G[N(x)]\rightarrow [K]$ be a proper coloring of $N(x)$ with at most $K$ colors. For every $x\in S$, choose a color $a_x\in [K]$ independently from the uniform distribution, and let
		$$Y=\{y\in Y_0: \forall x\in F_y, \chi_x(y)=a_x\}.$$
		Writing $q=K^{-k}$, we have $\mathbb{P}(y\in Y|y\in Y_0)=q$, so
		$$\mathbb{E}|Y|\geq \frac{1}{4}c_1qp^kM.$$
		If $y,y'\in Y$ such that $yy'$ is an edge, then $F_y$ and $F_{y'}$ are disjoint. Indeed, if $x\in F_{y}\cap F_{y'}$, then $\chi_x(y)=\chi_{x}(y')$, contradicting that $\chi_x$ is a proper coloring. From this, 
		$$\mathbb{P}(y,y'\in Y\mid y,y'\in Y_0, F_y\cap F_{y'}=\emptyset)=q^2.$$
		In particular, for every $yy'\in E(G)$, we can bound
		\begin{align*}
			\mathbb{P}(y,y'\in Y)&\leq q^2\sum_{\substack{F,F'\\ (y,F),(y',F')\in \cF \\ F\cap F'=\emptyset}} c_1^2\lambda_{y,F}\lambda_{y',F'}\mathbb{P}((y,F),(y',F')\text{ are excellent}) \\
			&\leq c_1^2q^2 p^{2k} \sum_{\substack{F,F'\\ (y,F),(y',F')\in \cF \\ F\cap F'=\emptyset}}\lambda_{y,F}\lambda_{y',F'}\\
			&\leq c_1^2q^2 p^{2k} \left(\sum_{F: (y,F)\in \cF}\lambda_{y,F}\right)\cdot \left(\sum_{F':(y',F')\in \cF}\lambda_{y',F'}\right)\\
			&\leq c_1^2q^2p^{2k}\Delta \sum_{F: (y,F)\in \cF}\lambda_{y,F}
		\end{align*}
		Therefore,
		\begin{align*}
			\mathbb{E}e(G[Y])&\leq \sum_{yy'\in E(G)}c_1^2q^2p^{2k}\Delta\sum_{F:(y,F)\in \cF}\lambda_{y,F}\\
			&\leq c_1^2q^2p^{2k}\Delta D\sum_{(y,F)\in \cF}\lambda_{y,F}\leq c_1^2q^2p^{2k}\Delta DM.
		\end{align*}
		Recall that $\Delta=C_1K^kD^{k-1}$, $q=K^{-k}$ and $p=c_0/D$, so 
		$$c_1^2q^2p^{2k}\Delta DM=c_0^{k}c_1^2C_1 q p^{k} M.$$
		Therefore, choosing $C_1=c_0^{-k}c_1^{-1}/8$, we have $\mathbb{E}e(G[Y])\leq \mathbb{E}|Y|/2$. Remove an endpoint of each edge in $G[Y]$, and let $I$ be the resulting set of vertices. Then $I\subset Y$ is an independent set of $G$ satisfying
		$$\mathbb{E}|I|\geq \frac{1}{2}\mathbb{E}|Y|\geq \frac{1}{8}c_1qp^kM.$$
		The main observation is that $G[S\cup I]$ is a $C_4$-free bipartite graph with vertex classes $S$ and $I$ such that every vertex in $I$ has degree exactly $k$. Indeed, every $y\in I\subset Y_0$ satisfies that $N(y)\cap S=F_y$, and every pair $y,y'\in I\subset Y_0$ satisfies that $|F_y\cap F_{y'}|\leq 1$, which ensures $C_4$-freeness. But then $\db(G)\leq k$ implies that 
		$$k|I|=e(G[I\cup S])\leq \frac{k(|I|+|S|)}{2},$$
		which gives that $|I|\leq |S|$. Hence,
		$$\frac{1}{8}c_1qp^kM\leq \mathbb{E}|I|\leq \mathbb{E}|S|\leq pn,$$
		from which 
		$$M\leq \frac{8n}{c_1qp^{k-1}}=\frac{8K^knD^{k-1}}{c_0^{k-1}c_1}=64c_0C_1K^knD^{k-1}=64c_0\Delta n\leq \frac{\Delta n}{8}.$$
		This finishes the proof.
	\end{proof}
	
	\begin{corollary}\label{cor:stable}
		Let $k,r\geq 2$ be integers, then there exist positive integers $C,t$ such that the following holds. Let $G$ be a nonempty $n$-vertex graph that is $K_r$-free, $(\beta,1/4)$-semiregular, and $\db(G)\leq k$. Assume that $G[N(x)]$ is $K$-colorable for every $x\in V(G)$. Then there exists $A\subset V(G)$ of size at least $3n/4$ such that $G[A]$ is $(\beta,1/2)$-semiregular, has maximum degree $D$, and $G[A]$ has the stable-$(CK^k\beta^{-t+1},D,t)$-property.
	\end{corollary}
	
	\begin{proof}
		Let $D_0$ be the maximum degree of $G$, and let $A\subset V(G)$ such that $|A|\geq 3n/4$ and $G[A]$ has the stable-$(CK^k,D_0,t)$-property. As $|A|\geq 3n/4$, we have that $G[A]$ is $(\beta,1/2)$-semiregular, and $G[A]$ has at least $\beta D_0n$ edges. But then $D\geq \beta D_0$, so by item 2 in \Cref{claim:stable_property}, $G[A]$ also has the stable-$(CK^k\beta^{-t+1},D,t)$-property.
	\end{proof}

    Our second main lemma is an analogue of \Cref{lemma:s-property}. We prove that a $(\beta,1/2)$-semiregular graph with the stable-$(L,D,s+1)$-property contains a large induced subgraph having the stable- $(L',D,s)$-property, where $L'$ is bounded by a function of $L$, $\beta$ and $s$.
	
	\begin{lemma}\label{lemma:stable-s-decrease}
		Let $s\in \mathbb{Z}^+$, then there exists $c_s>0$ such that the following holds. Let $L\geq 1$, and let $G$ be an $n$-vertex $(\beta,1/2)$-semiregular graph with  $\Delta(G)=D$ and the stable-$(L,D,s+1)$-property.
		Then there exists $X\subset V(G)$ such that $|X|\geq \Omega(\beta n/L)$, and $G[X]$ has the stable-$(c_s L^3/\beta^2,D,s)$-property.
	\end{lemma}
	
	\begin{proof}
		Let $\cH$ be a hypergraph on $V(G)$ witnessing that $G$ has the stable-$(L,D,s+1)$-property. First, we remove all vertices of $G$, whose degree in $\cH_2$ is more than
		$$\Delta:=4LD.$$
		As $|\cH_2|\leq LnD$, we removed at most 
		$$\frac{2LnD}{\Delta}=\frac{n}{2}$$
		vertices; let $U$ be the resulting set.
		For $x\in U$, let 
		$$Q_x=\{v\in U: \{x,v\}\in \cH_2\}$$
		be the neighbors of $x$ in $\cH_2$, then $|Q_x|\leq \Delta$.  We note that as $G$ is $(\beta,1/2)$-semiregular, we have $$e(G[U])\geq \beta Dn.$$
		
		Sample each vertex of $U$ with probability $p=1/(4D)$, and let $S$ be the set of sampled vertices. Moreover, further sample each element of $S$ with probability 
		$q=1/(4L),$ let $S^*$ be the resulting set. Let 
		$$U_0=S\cup \bigcup_{v\in S^*} Q_v,$$
		and let $X$ be the set of vertices $v\in U\setminus U_0$ such that $v$ has exactly one neighbor in $S$, and this one neighbor is in $S^*$. First, we show that $|X|$ is large with non-negligible probability. Let
		$$\lambda=\frac{\beta}{64L}.$$
		
		\begin{claim}
			We have
			$$\mathbb{P}\left(|X|\geq \lambda n\right)\geq \lambda.$$
		\end{claim}
		\begin{proof}
			For every $v\in U$, write $d(v)$ for the degree of $v$  in $G[U]$. Then $v\in X$ if and only if $Q_v\cap S^*=\emptyset$, $v\not\in S$, and $|N(v)\cap S^*|=|N(v)\cap S|=1$. Observe that as every edge in $\cH_2$ is an independent set in $G$, we have that $Q_v$ and $N(v)$ are disjoint. Thus, the events $$\{Q_v\cap S^*=\emptyset\text{ and }v\not\in S\}\quad\text{and}\quad\{|N(v)\cap S^*|=|N(v)\cap S|=1\}$$ are independent. The probability of the first event is $(1-p)(1-pq)^{|Q_v|}$, while the probability of the second event is $d(v)pq(1-p)^{d(v)-1}$. Therefore,
			$$\mathbb{P}(v\in X)=d(v)pq(1-p)^{d(v)}(1-pq)^{|Q_v|}\geq \frac{d(v)}{16DL}(1-p)^{D}(1-pq)^{\Delta}\geq 
			\frac{d(v)}{64DL}.$$
			From this,
			$$\mathbb{E}|X|=\sum_{v\in U}\mathbb{P}(v\in X)\geq \frac{1}{64DL}\sum_{v\in U}d(v)=\frac{1}{32DL}e(G[U])\geq \frac{\beta n}{32L}.$$
			As $|X|\leq n$, \Cref{lemma:simple_prob} applied with $c=1/2$ implies that 
			$$\mathbb{P}\left(|X|\geq \frac{\beta n}{64L}\right)\geq \frac{\beta}{64L}.$$
		\end{proof}
		
		Next, we define the hypergraph $\cH'$ on vertex set $X$ as follows. For $h=2,\dots,s$, let
		$$\cH'_h=\cH_h[X]\cup \{F\setminus \{x\}: F\in \cH_{h+1}\text{ and }\{x\}=F\setminus X=F\cap S^*\}.$$
		That is, we take all edges of $\cH[X]$ of size $h$, and in addition the truncations of those edges of size $h+1$ in $\cH$, that intersect $V(G)\setminus X$ in a single element of $S^*$. Every edge of $\cH'$ is an independent set in $G[X]$, as every edge of $\cH'$ is a subset of some edge of $\cH$. For every $F\in \cH_{h+1}$, the probability that some $x\in F$ is contained in $S^*$ is at most $(h+1)pq\leq (s+1)/D$. Thus, we have
		$$\mathbb{E}|\cH'_h|\leq |\cH_h|+\frac{s+1}{D}|\cH_{h+1}|\leq (s+2)LnD^{h-1}.$$
		By Markov's inequality, this gives
		$$\mathbb{P}\left(|\cH'_h|\geq \frac{s(s+2)}{\lambda}LnD^{h-1}\right)\leq \frac{\lambda}{s}.$$
		Therefore, with positive probability, the events
		$$|X|\geq \lambda n,\quad |\cH'_h|\leq \frac{s(s+2)}{\lambda}LnD^{h-1}\quad (2\leq h\leq s)$$
		hold simultaneously. Thus, we can fix an outcome of $S$ and $S^*$ such that all these events hold.
		
		Finally, we prove that for every $v\in X$, $\cH'[N(v)\cap X]$ has no independent set of size $s$. Recall that if $v\in X$, then there is a unique vertex $z$ such that $\{z\}= S^*\cap N(v)$. Let $J\subset N(v)\cap X$ of size $s$, and assume to the contrary that $J$ is an independent set in $\cH'$. As $\cH'$ contains all edges of $\cH[X]$ of size at most $s$, $J$ must be an independent set in $\cH$ as well. But as $J\cup \{z\}$ is not an independent set in $\cH[N(v)]$, there exists some edge $F$ in $\cH$ such that $F\subset J\cup\{z\}$. But then we must have $z\in F$. If $F$ has size 2, then $F=\{w,z\}\in \cH_2$ for some $w\in Q_z$. But this is impossible as $X$ is disjoint from $Q_z$, contradicting $w\in X$. On the other hand, if $|F|\geq 3$, then $F\setminus\{z\}$ is an edge of $\cH'$, and $F\setminus\{z\}\subset J$, contradicting that $J$ is an independent set. 
		
		Let $n'=|X|$, then $n'\geq \lambda n$, and for $h=2,\dots,s$,
		$$|\cH_h'|\leq \frac{s(s+2)}{\lambda^2}Ln'D^{h-1}=O_s(L^3/\beta^2) n'D^{h-1}.$$
		Thus, $G[X]$ has the stable-$(O_s(L^3/\beta^2),D,s)$-property. This finishes the proof.
	\end{proof}

    Our final lemma shows that if $G$ is a $K_r$-free degree-bounded graph in which the neighborhoods have small chromatic number, then $G$ contains a very large independent set. From this statement, \Cref{thm:main_color_precise} follows almost immediately. Indeed, proceeding by induction on $r$ ensures that the neighborhoods have small chromatic numbers, while repeatedly removing large independent sets provides a proper coloring of the graph with few colors.
	
	\begin{lemma}\label{lemma:independent_set}
		Let $k,r\geq 2$ be integers and $\tau>0$, then there exists $\eps>0$ such that the following holds if $n$ is sufficiently large. Let $G$ be an $n$-vertex $K_r$-free graph such that $\db(G)\leq k$. Assume that $G[N(x)]$ is $n^{\eps}$-colorable for every $x\in V(G)$. Then $G$ contains an independent set of size at least $n^{1-\tau}$.
	\end{lemma}
	
	\begin{proof}
		Let $t,C$ be the constants given by \Cref{cor:stable}, depending only on $k$ and $r$, and let $c_s$ be the constant given by \Cref{lemma:stable-s-decrease}, depending only on $s$. Choose $\eps>0$ sufficiently small with respect to $\tau$ and $t$, specified later. Furthermore, let $\rho=\rho(\eps,1/4)>0$ be the constant given by \Cref{lemma:semiregular_nonbip}. We assume that $n$ is sufficiently large with respect to $\eps,k,\tau$, which will be implicitly used in our calculations.
		
		Let $K=n^{\eps}$ and $\beta=1/K$, and define two sequences of positive real numbers $L_1,\dots,L_t$ and $x_2,\dots,x_t$ as follows. Let $L_t=CK^k\beta^{-t+1}=Cn^{(k+t-1)\eps}$, and for $s=t,t-1,\dots,2$, let 
		$$x_s=(\beta^{-t}L_t\dots L_s)^{2/\tau}\quad\text{and}\quad L_{s-1}=c_{s-1}L_s^3x_s^{s-2}\beta^{-2}/\rho.$$
		Clearly, we have $x_s,L_s\leq n^{O_{t,\tau}(\eps)}$.

		By \Cref{lemma:semiregular_nonbip}, we can choose a
		$(\beta,1/4)$-semiregular induced subgraph $G_t'$ of $G$
		on at least $\rho n$ vertices. We run the following algorithm, whose input is the graph $G_t'$.
		
		\begin{enumerate}
			\item \textbf{Initialize.}
			If $G_t'$ is edgeless, stop and output $G_t'$.
			Otherwise, use \Cref{cor:stable} to choose an induced subgraph $G_t$ on
			$n_t\ge |V(G_t')|/2$ vertices satisfying the conditions below.
			
			\item \textbf{Iterate.}
			For $s=t,t-1,\ldots,2$, assume that $G_s$ has
			$n_s$ vertices, $\Delta(G_s)=D_s$, is
			$(\beta,1/2)$-semiregular, and has the
			stable-$(L_s,D_s,s)$-property. Perform the following steps.
			
			\begin{enumerate}
				\item 
				By \Cref{lemma:stable-s-decrease}, we can choose $H=G_s[X]$
				such that $|X|=\Omega(\beta n_s/L_s)$ and $H$ has the
				stable-$(M_s,D_s,s-1)$-property, where
				\[
				M_s=\frac{c_{s-1}L_s^3}{\beta^2}.
				\]
				If $s=2$, stop and output $H$.
				
				\item 
				By \Cref{lemma:semiregular_nonbip}, choose a
				$(\beta,1/4)$-semiregular induced subgraph $G_{s-1}$
				of $H$ on $n_{s-1}\ge \rho|X|$ vertices.
				By item~3 of \Cref{claim:stable_property}, it has the
				stable-$(M_s/\rho,D_s,s-1)$-property.
				If $G_{s-1}$ is edgeless, stop and output it.
				Set $D_{s-1}=\Delta(G_{s-1})$.
				
				\item 
				If $x_sD_{s-1}<D_s$, set $G_t'\gets G_{s-1}$
				and \emph{restart} at step~1.
				
				Otherwise, $D_s/D_{s-1}\le x_s$, so item~2 of
				\Cref{claim:stable_property} gives the
				$$\text{stable-}(x_s^{s-2}M_s/\rho,D_{s-1},s-1)\text{-property}.$$
				By the definition of $L_{s-1}$, this is equal to the stable-$(L_{s-1},D_{s-1},s-1)$-property. Since
				$(\beta,1/4)$-semiregularity implies
				$(\beta,1/2)$-semiregularity, the required conditions are satisfied
				for the next iteration.
			\end{enumerate}
		\end{enumerate}
		
		We now analyze this algorithm. We first look at how the sizes and maximum degrees change within one \emph{cycle} of the algorithm, which refers to the steps executed before the algorithm is restarted.
		Observe that
		$$\frac{n_{t}}{n_{s-1}}\leq O_{t}\left((\rho\beta)^{-(t-s+1)}L_t\dots L_s\right)< \frac{1}{2}x_s^{\tau/2}.$$ 
		Therefore, if we restarted the cycle at some $s$, then we have 
		$$\frac{n_{t}}{n_{s-1}}\leq \frac{1}{2}\left(\frac{D_t}{D_{s-1}}\right)^{\tau/2}.$$
		In particular, this implies the following. Let $G_s^{(\ell)}$ denote the graph $G_s$ in the $\ell$-th cycle, and write $n_s^{(\ell)}=n_s$, $D_s^{(\ell)}=D_s$. Then
		$$\frac{n_t^{(\ell)}}{n_t^{(\ell+1)}}\leq \left(\frac{D_t^{(\ell)}}{D_t^{(\ell+1)}}\right)^{\tau/2},$$
		which implies by taking the telescopic product that
		$$\frac{n_{t}^{(1)}}{n_t^{(\ell)}}\leq \left(\frac{D_t^{(1)}}{D_t^{(\ell)}}\right)^{\tau/2}\leq n^{\tau/2}.$$
		Hence, as $n_t^{(1)}\geq \rho n/2$, each cycle starts with a graph of size at least $n^{1-2\tau/3}$. Moreover, as $n_t^{(\ell)}/n_s^{(\ell)}\leq n^{O(\eps)}$, by choosing $\eps$ sufficiently small, we also have $n_s^{(\ell)}\geq n^{1-3\tau/4}$ for every valid $s$ and $\ell$. But then the algorithm outputs some empty graph with at least $n^{1-\tau}$ vertices.
	\end{proof}
	
	Finally, we prove \Cref{thm:main_color_precise}.

	\begin{proof}[Proof of \Cref{thm:main_color_precise}]
		We proceed by induction on $r$. If $r=2$, the statement is trivial, so assume that $r\geq 3$ and that the statement holds for $r$ replaced by $r-1$. Let $f(N)$ denote the minimum independence number of a graph $H$ on $N$ vertices which is $K_r$-free and satisfies $\db(H)\leq k$. 
		
		Let $\tau=\beta/2$, and let $\eps$  be given by \Cref{lemma:independent_set}. Assume that $N$ is sufficiently large with respect to $\tau,\eps$, and fix an $N$-vertex $K_r$-free graph $H$ with $\db(H)\leq k$ and $\alpha(H)=f(N)$.  Observe that the neighborhood of every vertex in $H$ is $K_{r-1}$-free, so the neighborhood of every vertex can be colored with $N^{\eps}$ colors by our induction hypothesis. Hence, we can apply \Cref{lemma:independent_set} to $H$ to ensure that $f(N)\geq N^{1-\tau}$.
		
		Now color $G$ as follows. As long as $G$ has at least $n^{\beta/2}$ vertices,  remove one of its largest independent sets $I$, and color $I$ with a unique color. Finally, color each of the remaining vertices with a unique color. Assuming $n$ is sufficiently large, at each step we removed at least $|V(G)|^{1-\tau}$ vertices. Therefore, it is easy to see that this coloring uses at most $O_{\tau}(n^{\tau})+n^{\beta/2}\leq n^{\beta}$ colors, finishing the proof.  
	\end{proof}

    \section{Concluding remarks}

    In \Cref{thm:main_EH}, we proved that if $G=(A,B;E)$ is a bipartite graph with $|A|=|B|=n$, $\db(G)\leq k$ and $|E|\leq \frac{1}{2}|A||B|$, then there exists an anticomplete pair $(X,Y)$ with $X\subset A$, $Y\subset B$ and $|X|=|Y|\geq n^{1-o(1)}$. A careful inspection of our proof yields the more precise bound 
    $$|X|=|Y|\geq n\exp\left(-O_k((\log n)^{1-\frac{1}{2k}})\right).$$
   In contrast, in \Cref{thm:construction}, we constructed $n$-vertex bipartite graphs $G$ with $\db(G)=O(1)$ in which every equal sized anticomplete pair $(X,Y)$ across the vertex classes satisfies $|X|=|Y|\leq n(\log\log n)^{-1/2}$. Closing the gap between these bounds remains an interesting open problem.

   Similarly, \Cref{thm:main_color} shows that if $G$ is an $n$-vertex $K_r$-free graph with $\db(G)\leq k$, then $G$ can be properly colored with $n^{o(1)}$ colors. Our proof gives the more precise bound
   $$\chi(G)\leq \exp\left(O_{k,r}((\log n)^{1-b})\right),$$
   where $b=b(k,r)>0$ depends only on $k$ and $r$. In contrast, Pawlik et al. \cite{Pawlik} constructed $n$-vertex triangle-free segment intersection graphs  for arbitrarily large $n$ with chromatic number $\Omega(\log \log n)$. Since string graphs are degree-bounded, such graphs satisfy $\db(G)=O(1)$. Determining the correct order of growth between these lower and upper bounds is another natural problem.
   
   Several important degree-bounded families admit polylogarithmic upper bounds. Fox and Pach \cite{FP} proved that if $G$ is an $n$-vertex $K_r$-free string graph, we have $\chi(G)=(\log n)^{O(\log r)}$. String graphs form a subfamily of the graphs excluding all induced subdivisions of a suitable graph $H$. For every fixed $H$, Nguyen, Scott, and Seymour \cite{NSS2} proved that an $n$-vertex $K_r$-free graph $H$ excluding all induced subdivisions of $H$ has chromatic number $O_{r,H}((\log n)^{O(r)})$. Similarly, Bierbaumer and Buci\'c \cite{BB} proved that if $G$ is an $n$-vertex $K_r$-free graph avoiding a fixed tree $T$ as an induced subgraph, then $\chi(G)\leq O_{r,T}((\log n)^{O(r)})$. For this latter family, the celebrated Gy\'arf\'as-Sumner conjecture \cite{Gy1,Sumner} predicts the stronger bound $\chi(G)=O_{r,T}(1)$. 

    \section*{Acknowledgments}
\vspace{-10pt}
    IT acknowledges the support of the Swedish Research Council grant VR 2023-03375. We would like to thank Zach Hunter, Aleksa Milojevi\'c, and Benny Sudakov for many discussions about degree-bounded families, and for the construction presented in \Cref{lemma:simple_construction}. 
\vspace{-10pt}
    \section*{AI declaration}
    \vspace{-10pt}
    The proof of \Cref{thm:main_EH} is the result of interaction between the author and ChatGPT Astra. After finding a proof, the language model was able to adapt it very quickly to resolve \Cref{thm:main_color} as well.

\end{document}